\magnification=\magstep1
\input amstex
\documentstyle{amsppt}
\vsize=9truein
\def\Cirmin{\hbox{$\bigcirc\hbox{\kern-8pt\lower.5pt\hbox{{\tt -}}}$}}
\def\Cirtwo{\hbox{$\bigcirc\hbox{\kern-8pt\lower.5pt\hbox{{\tt 2}}}$}}
\def\Cirthree{\hbox{$\bigcirc\hbox{\kern-8pt\lower.5pt\hbox{{\tt 3}}}$}}
\def\Cirfour{\hbox{$\bigcirc\hbox{\kern-8pt\lower.5pt\hbox{{\tt 4}}}$}}
\def\Cirfive{\hbox{$\bigcirc\hbox{\kern-8pt\lower.5pt\hbox{{\tt 5}}}$}}
\define\sa{\operatorname{sa}}
\define\at{\operatorname{at}}
\define\rank{\operatorname{rank}}
\define\dist{\operatorname{dist}}

\NoBlackBoxes
\document

\topmatter
\title
{MASA's and Certain Type I Closed Faces of $C^*-$algebras}
\endtitle
\author
{Lawrence G. Brown}
\endauthor
\dedicatory
{Dedicated to the memory of George W. Mackey}
\enddedicatory
\footnote""{AMS subject classification:  46L05}

\abstract
{Let $A$ be a separable $C^*-$algebra and $A^{**}$ its enveloping $W^*-$algebra.  
A result of Akemann, Anderson, and Pedersen states that if $\{p_n\}$ is a
sequence of mutually orthogonal, minimal projections in $A^{**}$ such
that $\sum_k^\infty\ p_n$ is closed, $\forall k$,
then there is a MASA $B$ in $A$ such that each $\varphi_n|B$ is pure
and has a unique state extension to $A$, where $\varphi_n$ is the pure
state of $A$ supported by $p_n$.  We generalize this result in
two ways: We prove that $B$ can be required to contain an approximate
identity of $A$, and we show that the countable discrete space which
underlies the result cited can be replaced by a general totally disconnected
space.  We  consider two special kinds of type I closed faces, both related
to the above, atomic closed faces and closed faces with nearly closed
extreme boundary.  
One specific question is whether an atomic closed face
always has an ``isolated point''. We give a counterexample for this and also 
show that the answer is yes if the atomic face has nearly closed extreme
boundary.  We prove a complement to Glimm's theorem on type I $C^*-$algebras
which arises from the theory of type I closed faces.  One of our examples is a type I closed face which is isomorphic
to a closed face of every non-type I separable $C^*-$algebra and which is not
isomorphic to a closed face of any type I $C^*-$algebra.}
\endabstract
\endtopmatter
\document

\vfill\eject
\noindent
{\bf 0.  Introduction.}

This paper was inspired by the paper \cite{5} of C. Akemann, J. Anderson,
and G. Pedersen.  Much of the terminology used in this section is explained
in later sections.  To explain the connection with \cite{5}, we begin with:

\proclaim{Proposition 0.1}
Let $A$ be a $C^*-$algebra and $(p_n)$ a sequence of mutually
orthogonal, minimal (rank one) projections in $A^{**}$.  Let 
$p=\sum_1^\infty\ p_n$, and let $\varphi_n$ be the pure state supported by
$p_n$. Then either of the following hypotheses implies that $p$ is closed:
\roster
\item"(i)"
{\rm ([5, 2.7(1)}$\Rightarrow${\rm (2)]).} There is a strictly positive element $e$
such that each $\varphi_n$ is definite on $e$ and $\varphi_n(e)\to 0$.

\item"(ii)"
{\rm ([12, Lemma 3]).}  There is a strictly positive element $e$ such that
$\sum_1^\infty \varphi_n(e)<\infty$.
\endroster
\endproclaim

In circumstances similar to 0.1, \cite{5} proves the existence of a
MASA $B$ such that each $\varphi_n|_B$ has the unique extension property.
The hypotheses require that $A$ be non-unital.  It is known (see
[6, \S4]) that a non-unital $C^*-$algebra $A$ may have MASA's which do not
hereditarily generate $A$, or equivalently which do not contain an
approximate identity of $A$.  If the MASA constructed in \cite{5} does not
hereditarily generate $A$, the situation is intuitively unsatisfactory.
(See the first paragraph of [5, \S2].)

To investigate strengthening the result of \cite{5}, consider $\widetilde A$,
the result of adjoining an identity to $A$, and the pure state $\varphi_\infty$
defined by $\varphi_\infty(\lambda 1_{\widetilde A} +a)=\lambda$. The existence of a
MASA $B$ of $A$ such that each $\varphi_n|_B$, $1\le n<\infty$, has the
unique extension property and such that $B$ hereditarily generates $A$
is equivalent to the existence of a MASA $B_1$ of $\widetilde A$ such that each
$\varphi_n|_{B_1}$, $1\le n\le\infty$, has the unqiue extension property
$(B=B_1\cap A, B_1=\widetilde B)$.

Now the hypotheses of \cite{5} imply that $\sum_{n\in I} p_n$ is closed for
every subset $I$ of ${\Bbb N}$.  Thus $\{p_n:1\le n<\infty\}$ has properties
analogous to those of the discrete topological space ${\Bbb N}$.  But when
$p_\infty$, the support projection of $\varphi_\infty$, is added to the set,
the new set resembles the non-discrete space ${\Bbb N}\cup\{\infty\}$.  
Thus we seek
a generalization of the MASA result of \cite{5} based on a class of topological
spaces which includes ${\Bbb N}\cup\{\infty\}$.  We accomplish 
this in Corollary 2.4:
Let $A$ be a separable $C^*-$algebra and $X$ a totally disconnected,
second countable, locally compact Hausdorff space.  Assume that for each
$x$ in $X$, $p_x$ is a minimal projection in $A^{**}$, with associated
pure state $\varphi_x$, such that the $p_x$'s are mutually orthogonal and for
each closed (compact) subset $S$ of $X$, $\sum_{x\in S} p_x$ is the atomic part
of a closed (compact) projection, $p_S$, in $A^{**}$.  Then there is a MASA
$B$ of $A$ such that $B$ hereditarily generates $A$ and each $\varphi_x|_B$
has the unique extension property.  Moreover, each $p_S$ is in $B^{**}$.

If $X$ is general, the hypotheses of the above result may seem rather
stringent.  Partly in order to justify the generality of the result, we
attempt to investigate the circumstances in which the hypotheses will be
satisfied.  A first observation is that every element of $C_0(X)$ (respectively,
$C_b(X)$) gives rise to an element of $p_XA^{**}p_X$ which is strongly
$q$-continuous (respectively $q$-continuous) on $p_X$. (The concept ``$q$-continuous on $p$'' was defined in \cite{7}.  In \cite{11}, ``strongly 
$q$-continuous on $p$'' was defined and ``Tietze extension theorems'' for 
both kinds of relative $q$-continuity were given.)  Thus in Section 3 we give
some basic results and examples on the subject of how many relatively
$q$-continuous elements are supported by a given closed projection.

We also focus on a more specific question suggested by the theory of scattered
$C^*-$algebras:  Suppose that $p$ is an atomic closed projection in
$A^{**}$ and that $pA^*p$ is norm separable.  Is there a minimal projection
$p_0$ such that $p_0\le p$ and $p-p_0$ is closed?  Such a $p_0$ would give
an ``isolated point'' of the closed face $F(p)$ supported by $p$. This question
is related to the special case of 2.4 where $X$ is countable. Clearly, if
we seek to prove that certain conditions imply the hypotheses of 2.4, then we
must be able to prove that these conditions imply a positive answer to our
isolated point question.  Note also that when $X$ in 2.4 is countable,
then $p_X$ is atomic, and the words ``the atomic part of'' can be omitted.

We give a counterexample for the isolated point question in Section 3, but
we also give a positive result which has the following hypothesis (nearly
closed extreme boundary):
$$
[P(A)\cap F(p)]^-\subset\{0\}\cup[t,1]P(A)\ \text{ for some $t$ in }
(0,1].\tag NCEB
$$
Here $P(A)$ is the pure state space of $A$, and (NCEB) holds in particular
if the set of extreme points of $F(p)$ is (weak${}^*$) closed. Lest (NCEB)
seem unnatural or excessively strong, we point out a connection with [5,
\S4].   Circumstances not covered by 0.1
are actually considered in \cite{5}.
Suppose $\{\varphi_n:1\le n<\infty\}$ is a collection of mutually orthogonal
pure states such that $ \varphi_n \overset{w^*}\to\longrightarrow 0$ and each 
equivalence
class is finite.  With the additional assumption that there is a uniform
bound on the size of the equivalence classes, the authors of \cite{5} show
in \S4 that the needed conditions ([5, 2.7(2)]) are satisfied.  Without the
uniform boundedness hypothesis, we can show easily that [5, 2.7(2)] is
equivalent to (NCEB).

Our positive result, which is in Section 4, is roughly that if $p$ is a closed
projection satisfying (NCEB), then equivalence of pure states gives a proper
closed map from \newline
$[P(A)\cap F(p)]^-\setminus\{0\}$ onto a locally compact
Hausdorff space $X$.  If $pA^*p$ is norm separable, then $X$ is countable
and hence scattered. In general, $X$ need not even be totally disconnected,
of course.

There are some technicalities involving direct integral theory required in
order to prove that closed subsets of $X$ give rise to closed projections.
This is what leads us to the study in Section 5 of type I closed faces, where
the face $F(p)$ is called type I if $pA^{**}p$ is a type I $W^*-$algebra.
Obviously every atomic face is type I, and also $F(p)$ is type I when $p$ (is
closed and) satisfies (NCEB), at least if $A$ is separable. Our results on
type I closed faces are only rudimentary, and we think the concept is worthy
of further study.

Partly because theorems are not always discovered in logical order, our
efforts to expand on the results of \cite{5} have led us in several directions.
The different parts of this paper, though closely related, do not mesh
perfectly.  In Section 7 we attempt to exhibit the formal relationships among
the previous sections.  The earlier sections can in large part be read
independently of one another, except that Section 6 is a continuation of
Section 4 relying on Section 5.  The promised complement to Glimm's theorem
is Proposition 5.11.

A preliminary preprint of this paper was circulated several years ago.  Some
results overlapping with Section 3 have been independently found by
E. Kirchberg (cf\. [22, Lemma 2.3]).

\noindent
{\bf 1.  Preliminaries}.

$A$ will always denote a $C^*-$algebra and $A^{**}$  its enveloping
$W^*-$algebra.  For $h$ in $A^{**}_{\sa}$ and $F$ a Borel set in
$\Bbb R$ $E_F(h)$ denotes the spectral projection of $h$ for $F$.  
For many of our proofs $A$ must be separable, but we rarely
require that  $A$ be unital.  $S(A)$ is the state space of $A$,
$P(A)$ the pure state space, and $Q(A)$ the quasi-state space (the
set of positive functionals of norm at most $1$).  If $p$ is a projection in
$A^{**}$, $F(p)=\{\varphi\in Q(A):\varphi(1-p)=0\}$, the norm closed face
of $Q(A)$ supported by $p$.  (Elements of $A^*$ are regarded as functionals
on $A^{**}$ without notice.) Topological terminolgy regarding $A^*$ refers
to the weak${}^*$ topology unless the contrary is explicitly indicated.  A projection $p$
in $A^{**}$ is called {\it open} (\cite{1}) if it is the support projection of a
hereditary $C^*-$subalgebra of $A$ and {\it closed} if $1-p$ is open. Effros
proved in \cite{{15}, Theorem 4.8} (cf\. [25, 3.10.7]) that $p$
is closed if and only if $F(p)$ is closed, and if so we follow
the usual abuse of notation and call $F(p)$ a ``closed face of $A$''.
Also, $p$ is called {\it compact} (\cite{4}) if $F(p)\cap S(A)$ is closed or
equivalently if $p$ is closed in $\widetilde A^{**}$, where $\widetilde A$
is the result of adjoining an identity to $A$.

The {\it reduced atomic representation}, $\pi$, of $A$ is $\oplus_i \pi_i$,
where $\{\pi_i\}$ contains one representative of each unitary equivalence
class of irreducible representations. Denote by $z_{\at}$ the central
projection in $A^{**}$ that supports $\pi$. 
Thus  $z_{\at}A^{**}\cong\oplus_iB(H_{\pi_i})$,
and $(1-z_{\at})A^{**}$ has no type I factor direct summands.  The
{\it atomic part} of an element $x$ of $A^{**}$ is $z_{\at}x$, $x$ is
{\it atomic} if $x=z_{\at}x$, $F(p)$ is {\it atomic} if $p$ is atomic, etc.
Also pure states $\varphi$ and $\psi$ are called {\it equivalent} if
the irreducible representations $\pi_\varphi$, $\pi_\psi$ are unitarily
equivalent.

We want to comment further on Proposition 0.1.  In fact
2.7(1) of \cite{5} actually states that $\sum \varphi_n(e)<\infty$ rather
than $\varphi_n(e)\to 0$.  However, the proof in \cite{5} that (1) implies
(2) uses only the weaker hypothesis, so that it is correct to attribute
0.1(i) to \cite{5}.  (Unfortunately, when he was writing \cite{12}, the
author had not yet read the proofs in \cite{5}.)  Here is a generalization
of 0.1:

\proclaim{Lemma 1.1}
Let $(p_n)$ be a sequence of mutually orthogonal minimal projections in
$A^{**}$ and $p=\sum_1^\infty p_n$.  If, $\forall a\in A$, $\pi^{**}(p)
\pi(a)\pi^{**}(p)$ is a compact operator on $H_\pi$, where $\pi$ is the
reduced atomic representation of $A$, then $p$ is closed.
\endproclaim

The proof of 1.1 and the fact that it implies 0.1(ii) is identical to the
proof of Lemma 3 in \cite{12}. Lemma 1.1 implies 0.1(i) because in that case
$\pi^{**}(p)\pi(e)\pi^{**}(p)$ is a diagonal operator whose matrix
elements approach zero.  (If $A$ is $\sigma$-unital, it is enough to
verify the compactness for a strictly positive element of $A$, as shown
in \cite{12}.)  Lemma 1.1 also applies under the Standing Assumptions of
[5, \S4], since then $\pi^{**}(p)\pi(a)\pi^{**}(p)$ is a block-diagonal
operator with bounded block size - in particular it is a $(2N+1)$-diagonal
operator.

Despite this, we offer the following new proof of 0.1(i), which may be
instructive:

The hypothesis that $\varphi_n$ is definite on $e$ is equivalent to $p_ne=ep_n$. Thus if $\lambda_n=\varphi_n(e)$, then $p_n\le E_{\{\lambda_n\}}(e)$.
Let $\epsilon_k=\sup\{\lambda_n:n>k\}$ and
$q_k=\sum_1^k p_n \vee  E_{[0,\epsilon_k]}(e)$.  Since $[0,\epsilon_k]$ is a closed set,  $E_{[0,\epsilon_k]}(e)$ is
closed, and thus 
\cite{{1}, Theorem II.7} implies that $q_k$ is closed.  Since $E_{\{0\}}(e)=0$,
$p=\wedge_1^\infty q_k$, and \cite{{1}, Proposition II.5} implies $p$ is
closed.

\proclaim{Theorem 1.2}
Let $(p_n)$ be a sequence of mutually orthogonal minimal projections in
$A^{**}$ and $p=\sum_1^\infty p_n$.  Then the following are equivalent:
\roster
\item"(i)"
Every subprojection of $p$ in $A^{**}$ is closed.
\item"(ii)"
$\sum_{n\in I} p_n$ is closed for each subset $I$ of ${\Bbb N}$.
\item"(iii)"
$\sum_k^\infty p_n$ is closed, $\forall k$ {\rm (cf. [5, 2.7(2)])}.
\item"(iv)"
$\pi^{**}(p)\pi(A)\pi^{**}(p)\subset\Cal K(H_\pi)$, where $\pi$ is the
reduced atomic representation of $A$.
\endroster
\endproclaim

\demo{Proof}
(i)$\Rightarrow$(ii)$\Rightarrow$(iii) is obvious.

(iii)$\Rightarrow$(iv):  Let $q_k=\sum_k^\infty p_n$.  Then $F(q_k)$ is a
closed subset of $Q(A)$ and $\underset 1\to{\overset\infty\to\bigcap} F(q_k)=
\{0\}$.  By this and [5, 2.3], $(q_k)$ approaches infinity in the sense
of [5].  By
definition, $\Vert aq_k\Vert\to 0$, $\forall a\in A$.  Therefore
$\pi(a)\pi^{**}(p-q_k)\to\pi(a)\pi^{**}(p)$ in norm. Since
$\pi^{**}(p-q_k)$ is a finite rank operator, this implies
$\pi(a)\pi^{**}(p)$, and {\it a fortiori} $\pi^{**}(p)\pi(a)\pi^{**}(p)$, is
compact.

(iv)$\Rightarrow$(i):  Assume $p'\in A^{**}$ and $p'\le p$.  Then
$\pi^{**}(p')\pi(A)\pi^{**}(p')\subset\Cal K(H_\pi)$.  Clearly there
are mutually orthogonal minimal projections $p'_n$ such that $p'=\sum_1^\infty
p_n'$.  Thus $p'$ is closed by 1.1.
\enddemo

Perhaps it should also be mentioned that if $I$ in 1.2(ii) is finite, then
$\sum_{n\in I} p_n$ is finite rank and hence compact ([1, Corollary II.8]).

If $\varphi$ is in $P(A)$ and $B$ is a $C^*-$subalgebra of $A$, we say
that $\varphi|_B$ has the {\it unique extension property} (UEP) if
$\varphi|_B\in P(B)$ and $\varphi$ is the only element of $S(A)$ which extends
$\varphi|_B$.  The next proposition is probably not original (see
[5, p. 267]).

\proclaim{Proposition 1.3}
Assume $p$ is  minimal projection in $A^{**}$ and $\varphi$ is the associated
pure state.  If $B$ is a $C^*-$subalgebra of $A$, then $\varphi|_B$ has
(UEP) if and only if $p$ is in $B^{**}$.
\endproclaim

\demo{Proof}
Of course $B^{**}$ is identified with the weak* closure of $B$ in $A^{**}$.
First assume (UEP) and let $q$ be the support projection of $\varphi|_B$,
so that $q$ is a minimal projection in $B^{**}$. If $\psi$ is in $S(A)$ and
$\psi(q)=1$, then $\psi|_B$ is in $F(q)\cap S(B)$, and hence
$\psi|_B=\varphi|_B$.  By (UEP), $\psi=\varphi$.  Thus we have shown that
$F(q)$, computed in $A^*$, is one dimensional, and this clearly implies
$q= p$.

Conversely, assume $p\in B^{**}$.  Since $p$ is minimal in $A^{**}$, it is 
clearly minimal in $B^{**}$.  Since $\varphi|_B(p)=1$ and
$\Vert\varphi|_B\Vert\le\Vert\varphi\Vert=1$, $\varphi|_B$ is a state
supported by $p$.  Therefore $\varphi|_B\in P(B)$.  If $\psi\in S(A)$ and
$\psi|_B=\varphi|_B$, then $\psi$ and $\varphi$ agree also on $B^{**}$.
Thus $\psi(p)=\varphi(p)=1$, $\psi$ is supported by $p$, and hence $\psi=\varphi$.
\enddemo

Recall the condition (NCEB), which was defined in Section 0 for any
projection $p$ in $A^{**}$.  It is also convenient to have a name for the
special case of (NCEB) where $t=1$.
$$
[P(A)\cap F(p)]^-\subset \{0\}\cup P(A). \tag CEB
$$
The phrase ``closed extreme boundary'', is accurate only when $p$ is closed,
but the main uses of (NCEB) and (CEB) are for projections known
{\it a priori} to be closed.

\proclaim{Theorem 1.4}
Let $(p_n)$ be a sequence of mutually orthogonal minimal projections in
$A^{**}$, $(\varphi_n)$ the associated sequence of pure states, and
$p=\sum_1^\infty p_n$.  If the equivalence classes of $\{\varphi_n\}$ are
finite and $\varphi_n\overset{w^*}\to\longrightarrow 0$, then the
following are equivalent:
\roster
\item"(i)"
$\quad\sum_k^\infty\ p_n$ is closed, $\forall k$.
\item"(ii)"
$\quad p$ satisfies (NCEB).
\item"(iii)"
$\quad p$ satisfies (CEB).
\item"(iv)"
$\quad [F(p)\cap P(A)]^-\subset\{0\}\cup[t,1]S(A)$ for some $t$
in $(0,1]$.
\endroster
\endproclaim

\demo{Proof}
Let $\Gamma_1,\Gamma_2,\dots$ be the equivalence classes of $\{\varphi_n\}$,
and let $q_i=\sum_{\varphi_n\in\Gamma_i} p_n$. Thus each $q_i$ is a
finite rank, and hence compact, projection in $A^{**}$.

(iv)$\Rightarrow$ (i): For this it is clearly permissible to simplify
the notation by assuming $k=1$.  Thus we need to show that $p$, which is
$\sum_i q_i$, is closed.  According to Proposition 4.2 of \cite{5}, for
this it is sufficient to show that $(q_i)$ approaches infinity.  Let $U$ be
a convex neighborhood of $0$ in $A^*$. We need to find $i_0$ such that
$F(q_i)\subset U$ for $i\ge i_0$. By the Krein-Milman theorem, it is
sufficient to show  $F(q_i)\cap P(A)\subset U$ for $i\ge i_0$. If this is
false we can find nets $(\psi_j)$ and $(i_j)$ such that 
$\psi_j\in F(q_{i_j})\cap
P(A)$, $i_j\to\infty$, $\psi_j\to\psi$, and $\psi\not=0$.

Let $\pi$ be the reduced atomic representation of $A$ and $H_j$ the range
of $\pi^{**}(q_{i_j})$.  Thus $\dim\ H_j=|\Gamma_j|$.  If $\dim\ H_j=1$ for
arbitrarily large $j$, then $\psi_j=\varphi_n$ for $\varphi_n\in\Gamma_{i_j}$;  
and we already know $\varphi_n\to 0$.  Thus we may assume
$\dim\ H_j\ge 2$, $\forall j$.  Then we can find unit vectors
$u_j$, $v_j'$, $v_j''$ in $H_j$ such that $(v_j',v_j'')=0$, the pure
states $(\pi(\cdot)v_j',v_j')$ and $(\pi(\cdot)v_j'',v_j'')$ are in
$\Gamma_{i_j}$, and $\psi_j=(\pi(\cdot)u_j,u_j)$.  Choose a unit vector $v_j$
in $\text{span}\{v_j',v_j''\}$ such that $(v_j,u_j)=0$, and let
$\theta_j=(\pi(\cdot)v_j,v_j)$.   Then $\theta_j\to 0$.  This follows from
Lemma 4.1 of \cite{5}, with the $N$ of \cite{5} being $2$, or it can be
proved directly using an argument similar to the one below.  Let
$f_j=(\pi(\cdot)u_j,v_j)$, which is an element of $A^*$.  By the Schwarz inequality,
$|f_j(a)|\le\Vert\pi(a^*)v_j\Vert=\theta_j(aa^*)^{1/2}$.  Therefore
$f_j\to 0$. Then if $w_j=ru_j+sv_j$, with $|r|^2+|s|^2=1$, and
$\rho_j=(\pi(\cdot)w_j,w_j)$, we see that $\rho_j\in F(p)\cap P(A)$, and
$\rho_j\to|r|^2\psi$.  We can choose $r,s$ such that $0<|r|^2\Vert\psi\Vert<t$,
in contradiction to (iv).

(i)$\Rightarrow$(iii):  Suppose $\psi_j\in F(p)\cap P(A)$ and $\psi_j\to\psi$.
Then for each $j$ there is $i_j$ such that $\psi_j\in F(q_{i_j})$.  If
$i_j\nrightarrow\infty$, then by passing to a subnet we may assume $i_j=i$, $\forall j$.
Then it is easy to see that $\psi\in F(q_i)\cap P(A)$.  (Each $\psi_j$
is a vector state coming from $H_j$ and the unit sphere of $H_j$ is norm
compact.)  If $i_j\to\infty$, then for each $k$, $\psi_j\in F(\sum_k^\infty 
p_n)$ for sufficiently large $j$.  By (i), $\psi\in F(\sum_k^\infty p_n)$. 
Since $\underset{k=1}\to{\overset\infty\to\bigwedge}\sum_k^\infty\ p_n=0$, $\psi=0$.

(iii)$\Rightarrow$(ii)$\Rightarrow$(iv) is obvious.
\enddemo

\noindent
{\bf 2.  Existence of MASA's}

\proclaim{Lemma 2.1}
Let $A$ be a $C^*-$algebra and $\widetilde A$ the result of adjoining a new
identity to $A$ (i.e., $\widetilde A\cong A\oplus{\Bbb C}$ if $A$ is already
unital).  Let $\varphi_\infty$ in $P(\widetilde A)$ be defined by
$\varphi_\infty(\lambda 1_{\widetilde A}+a)=\lambda$.   Assume $B_1$ is a unital $C^*-$subalgebra
of $\widetilde A$ such that $\varphi_{\infty}|_{ B_1}$ has $(UEP)$, and let
$B=B_1\cap A$.  Then $B$ hereditarily generates $A$ and $B^{**}=
B_1^{**}\cap A^{**}$.
\endproclaim

\demo{Proof}
That $B^{**}=B_1^{**}\cap A^{**}$ follows, for example, from
general Banach space theory and the fact that $B_1/B$ is finite dimensional.
Now let $p_\infty$ be the support projection of $\varphi_\infty$.  Then
$p_\infty\in B_1^{**}$ by 1.3.  Since $1_{\widetilde A}\in B_1^{**}$, $1_{\widetilde A}-p_\infty$ is also in
$B_1^{**}$, and of course $1_{\widetilde A}-p_\infty$ is the identity of $A^{**}$.  Thus
$1_{\widetilde A}-p_\infty\in B^{**}$, and this implies that $B$ hereditarily
generates $A$.
\enddemo

\proclaim{Lemma 2.2}
Let $A$ be a separable unital $C^*-$algebra, $p$ a closed projection in
$A^{**}$, and $q$ an open projection in $A^{**}$ such that $q\ge p$.  Let
$B=\text{her}(q)$, the hereditary $C^*-$subalgebra of $A$ supported by
$q$, and let $U$ be a neighborhood of $F(p)\cap S(A)$ in $S(A)$.  Then there
is a closed projection $p'$ in $B^{**}$ such that $p'p=0$ and
$\varphi(p')=0$ implies $\varphi\in U$ for $\varphi$ in $S(B)$.
\endproclaim

\demo{Proof}
As usual, we identify $B^{**}$ with $qA^{**}q$ and $S(B)$ with $\{\varphi
\in S(A):\varphi(q)=1\}$.  (The weak${}^*$ topologies of $A^*$ and $B^*$ agree
on $S(B)$.)  By Akemann's Urysohn lemma, \cite{{2}, Theorem 1.1}, there is $a$ in $A_{\sa}$ such
that $p\le a\le q$.  Then $a\in B$.  Let $C=\text{her}(q-p)$, and let
$e$ be a strictly positive element of $C$.  Let $b=a-aea$.  Then, by an
argument of Akemann [3, 1.1], $E_{\{1\}}(b)=p$.  Let $p_n'=E_{(-\infty,1-n^{-1}]}
(b)$, where the spectral projection is computed in $B^{**}$.  We claim that
for $n$ sufficiently large the choice $p'=p_n'$ suffices.  If not, for
each $n$ there is $\varphi_n$ in $S(B)$ such that $\varphi_n(p_n')=0$ and
$\varphi_n\notin U$. Then $\varphi_n$ is supported by
$q-p_n'=E_{(1-n^{-1},1]}(b)\le E_{[1-n^{-1},1]}(b)$.  Let $\varphi$ be a
cluster point of $(\varphi_n)$ in $S(A)$.  Then since each
$E_{[1-n^{-1},1]}(b)$ is closed in $A^{**}$, $\varphi$ is supported by
$\underset{n=1}\to{\overset\infty\to\bigwedge} E_{[1-n^{-1},1]}(b)=
E_{\{1\}}(b)=p$.  Therefore $\varphi\in F(p)\cap S(A)$, a contradiction
since $\varphi_n\notin U$.
\enddemo

\proclaim{Theorem 2.3}
Let $A$ be a separable $C^*-$algebra and $X$ a second countable, totally
disconnected, locally compact Hausdorff space.  Assume that for each $x$ in
$X$, $p_x$ is an atomic projection in $A^{**}$, the $p_x$'s are mutually
orthogonal, and for every closed (compact) subset $S$ of $X$ there is a
closed (compact) projection $p_S$ such that $z_{\at}p_S=\sum_{x\in S} p_x$.
Then there is a MASA $B$ in $A$ such that $B$ hereditarily generates $A$
and each $p_S$ is in $B^{**}$.
\endproclaim

\demo{Proof}
First we reduce to the case $A$ unital, $X$ compact. To do this, let
$\widetilde A$ be the result of adjoining a new identity to $A$, and let
$\widetilde X=X\cup\{\infty\}$ be the one point compactification.  If we let
$p_\infty$ be as in 2.1, all hypotheses of the theorem are satisfied for
$\widetilde X$, $\widetilde A$.  (If  $S$  is a compact subset of $X$, then
$p_S$ is compact in $A^{**}$ and hence closed in $\widetilde A^{**}$.  Any
other closed subset of $\widetilde X$ is $S\cup\{\infty\}$ for some closed
subset $S$ of $X$.  The fact that $p_S$ is closed in $A^{**}$ implies
that $p_S+p_\infty$ is closed in $\widetilde A^{**}$.  Since $\widetilde A$
is unital, ``closed''  and ``compact'' mean the same for projections in
$\widetilde A^{**}$.)  If $B_1$ satisfies the conclusion of the theorem
for $\widetilde A$, $\widetilde X$, then by 2.1, $B_1\cap A$ satisfies the
conclusions of the theorem for $A$, $X$.  Thus from now on we assume $A$
unital and $X$ compact.

Let $C$ be the usual middle-thirds Cantor set in $[0,1]$.  Then there is a
one-to-one continuous function $f:X\to C$.  We will let $\alpha$ and $\beta$ 
denote finite strings of $+$'s and $-$'s, and  $|\alpha|$ denote
the length of $\alpha$.  Let $C_+$, $C_-$ be the right and left halves of
$C$, $C_{++}$, $C_{+-}$ the right and left halves of $C_+$, etc.
Let $p_\alpha=p_{f^{-1}(C_\alpha)}$, and let $F_\alpha=S(A)\cap F(p_\alpha)$.
Note that $p_{\alpha+}p_{\alpha-}=0$ and $p_\alpha=p_{\alpha+}+
p_{\alpha-}$.  This follows, for example, from the theory of universally
measurable elements of $A^{**}$, [25, 4.3] and the fact that the relations
are satisfied by the atomic parts of the projections.  Let $e$ be a strictly
positive element of $\text{her}(1-p_X)$.  
We are going to construct recursively
$b_\alpha$ in $A_+$ and an open projection $q_\alpha$ in $A^{**}$ such that:
\roster
\item"1."
$b_{\alpha+}b_{\alpha-}=0$

\item"2."
$p_\alpha\le q_\alpha\le E_{\{1\}}(b_\alpha)$

\item"3."
$b_{\alpha+},b_{\alpha-}\in \text{her}(q_\alpha)$.  (Thus $b_\alpha
b_{\alpha\pm}=b_{\alpha\pm}$.)

\item"4."
If $\varphi$ in $S(A)$ is supported by $E_{\{1\}}(b_{\alpha\pm})$,
then $\varphi(e)< |\alpha|^{-1}2^{-|\alpha|}$.
\endroster

Fix non-negative functions $g_+$, $g_-$ in $C([-1,1])$ such that $g_+=1$
on $\left[\frac23, 1\right]$, $g_+$ is supported on $\left[\frac13, 1\right]$,
$g_-=1$ on $\left[-1,-\frac23\right]$, and $g_-$ is supported on
$\left[-1,-\frac13\right]$.

Step 1, $|\alpha|=1$.  Then 3 and 4 are vacuous.  Choose $a$ in $A_{\sa}$ such that
$-1\le a\le 1$, $p_-\le E_{\{-1\}}(a)$, and $p_+\le E_{\{1\}}(a)$.  This
is easily accomplished by \cite{{2}, Theorem 1.1} and the continuous functional calculus.
Let $b_\pm =g_\pm(a)$, $q_+=E_{(\frac23,1]}(a)$, and
$q_-=E_{[-1,-\frac23)}(a)$.

Step $k$, $|\alpha|=k>1$. We construct $b_{\beta\pm}$, $q_{\beta\pm}$ for each
$\beta$ with $|\beta|=k-1$, assuming of course that $b_\beta$, $q_\beta$
have already been constructed.  Apply 2.2 to find a closed projection
$p'$ in $\text{her}(q_\beta)^{**}$ such that $p'p_\beta=0$ and if
$\varphi$ in $S(A)$ is supported by $q_\beta$ and $\varphi(p')=0$, then
$\varphi(e)<|\beta|^{-1}2^{-|\beta|}$.  Next choose  $a$ in
$\text{her}(q_\beta)$ such that $-1\le a\le 1$, $p'\le E_{\{0\}}(a)$, and
$p_{\beta\pm}\le E_{\{\pm 1\}}(a)$.  The existence of $a$ could be deduced
from [11, 3.43], but it is more elementary to apply Akemann's Urysohn
lemma for $\text{her}(q_\beta)$ twice to obtain $a_1$ and $a_2$ with
$p_{\beta+}\le a_1\le 1-(p'+p_{\beta-})$ and $p_{\beta-}\le a_2\le
1-(p'+p_{\beta+})$.  Then let $a=a_1-a_2$.  Then let $q_{\beta+}=
E_{(\frac23,1]}(a)$, $q_{\beta-}=E_{[-1,-\frac23)}(a)$, and
$b_{\beta\pm}=g_\pm(a)$.

Now $\{b_\alpha\}$ is commutative, since for $\alpha\not=\alpha'$ either
$b_\alpha b_{\alpha'}=0$, $b_\alpha b_{\alpha'}=b_{\alpha'}$, or 
$b_\alpha b_{\alpha'}=b_\alpha$.  Let $B$ be any MASA containing all 
$b_\alpha$'s.
If $p_\alpha'=E_{\{1\}}(b_\alpha)$, then $p_\alpha'\in B^{**}$.  Note
that $p_{\alpha_1}'p_{\alpha_2}'=0$ if $|\alpha_1|=|\alpha_2|$ and
$\alpha_1\not=\alpha_2$ and that $p_\alpha'\ge p_\alpha$.

We show that $p_X\in B^{**}$ by proving $p_X=\underset{n=1}\to{\overset
\infty\to\bigwedge} \underset{|\alpha|=n}\to\bigvee p_\alpha'$.  Clearly
the latter is at least $p_X$.  Suppose $\varphi\in S(A)\cap F(
\underset{|\alpha|=n}\to\bigvee p_\alpha')$.  Let $\varphi_\alpha=
p_\alpha'\varphi p_\alpha'$.  Then $\underset{|\alpha|=n}\to\sum \Vert\varphi_
\alpha\Vert=1$, $\varphi_\alpha(e)<(n-1)^{-1}2^{-(n-1)}
\Vert\varphi_\alpha\Vert$, by 4, and $\varphi\le 2^n \underset{|\alpha|=n}\to
\sum \varphi_\alpha$.  Therefore $\varphi(e)< 2(n-1)^{-1}$.  If the above is
true for all $n$, then $\varphi(e)=0$ and hence $\varphi\in F(p_X)$.

Finally, to show that every $p_S$ is in $B^{**}$, note that every closed
subset of $C$ is the intersection of a sequence of clopen sets and every
clopen set is the union of finitely many $C_\alpha$'s.  Thus it is sufficient
to show that each $p_\alpha$ is in $B^{**}$.  We do this by showing that
$p_\alpha=p_X\wedge p_\alpha'$.  This follows from $p_\alpha' \ge p_\alpha$,
$p_\alpha' p_\beta = 0$ if $|\alpha| =|\beta|$ and $\alpha \ne \beta$,
and $p_X = \sum_{|\beta|=|\alpha|} p_\beta$.
\enddemo

\proclaim{Corollary 2.4}
Assume the hypotheses of 2.3 and also that each $p_x$ is a minimal projection
in $A^{**}$.  Let $\varphi_x$ be the pure state supported by $p_x$.  Then
if $B$ is the MASA of 2.3, $\varphi_x|B$ has the unique extension property,
$\forall x\in X$.
\endproclaim

\demo{Proof}
Combine 2.3 and 1.3, and note that $p_x=p_S$ for $S=\{x\}$.
\enddemo

\example{Remark 2.5} Since the construction of the MASA in 2.3 requires 
only the $p_S$'s, we could start with a more general, but also more abstract,
setup, an assignment $S \mapsto p_S$, for $S$ closed, such that:

(i) $p_{S_1}p_{S_2} =p_{S_2}p_{S_1}$,

(ii) $p_{\emptyset} =0$,

(iii) $p_{S_1 \cup S_2} = p_{S_1} \vee p_{S_2}$,

(iv) $p_{\cap_1^{\infty} S_n} = \bigwedge_1^{\infty} p_{S_n}$, and

(v) $p_S$ is closed and $S$ compact implies $p_S$ compact.

\noindent Because of our assumption that $X$ is totally disconnected,
condition (i) is redundant.  These conditions do not imply that
$z_{\at}p_S = \sum_{x\in S}z_{\at}p_{\{x\}}$, and this last property
is not needed to construct the MASA.  It was used in the proof of 2.3 to
prove conditions (iii) and (iv).

Another alternative formulation, using relative $q-$continuity, appears
below in 7.1 (see also 7.5).  The hypotheses actually used in 2.3 and 2.4 imply a stronger
relationship between the structure of $F(p_X)$ and the space $X$.
\endexample

\noindent
{\bf 3.  Relative $q$-continuity}

Let $p$ be a closed projection in $A^{**}$ and $h$ an element of
$pA_{\sa}^{**}p$.  Then $h$ is called  {\it $q$-continuous on $p$} (\cite{7}) if
$E_F(h)$ is closed for every closed subset $F$ of ${\Bbb R}$, where the
spectral projection is computed in $pA^{**}p$, and  $h$ is called
{\it strongly $q$-continuous on $p$} (\cite{11}) if in addition,
$E_F(h)$ is compact when $F$ is closed and $0\notin F$.  It was shown in
[11, 3.43] that $h$ is strongly $q$-continuous on $p$ if and only if
$h=pa$ for some $a$ in $A_{\sa}$ such that $pa=ap$, and if $A$ is
$\sigma$-unital, then $h$ is $q$-continuous on $p$ if and only if
$h=px$ for some $x$ in $M(A)_{\sa}$ such that $px=xp$.

It was neglected in \cite{11} to give any serious examples
or discussion of how extensive is the set of relatively $q$-continuous
elements.  For general $h$ in $pA^{**}p$ let us say that $h$ is
$q$-continuous or strongly $q$-continuous on $p$ if both Re $h$ and
Im $h$ are.  Let $SQC(p)=\{h\in pA^{**}p:h \text{ is strongly $q$-continuous
on } p\}$, and let $QC(p)=\{h\in pA^{**}p:h \text{ is $q$-continuous on }
p\}$.  By [11, 3.45], $SQC(p)$ is a $C^*-$algebra, and if $A$ is
$\sigma$-unital, $QC(p)$ is also a $C^*-$algebra.  We say that $p$ satisfies
$(MSQC)$ (many strongly $q$-continuous elements) if $SQC(p)$ is
$\sigma$-weakly dense in $pA^{**}p$ and $p$ satisfies $(MQC)$ if $QC(p)$ is
$\sigma$-weakly dense in $pA^{**}p$.  The von Neumann and Kaplansky density
theorems give many equivalent formulations of $(MSQC)$, and also
$(MQC)$ if $A$ is $\sigma$-unital.  As for the other extreme, we always
have ${\Bbb C}p\subset QC(p)$ and $0\in SQC(p)$.  We will show that
$QC(p)$ and $SQC(p)$ need not be any bigger.  Of course, $QC(p)=SQC(p)$ if
and only if $p$ is compact.

\proclaim{Theorem 3.1}  If $p$ is a closed projection in $A^{**}$, then the
following are equivalent:
\roster
\item"1."
$p$ satisfies $(MSQC)$.

\item"2."
$pAp=SQC(p)$.

\item"3."
$pAp$ is an algebra.

\item"4."
$pAp$ is a Jordan algebra.

\item"5."
$F(p)$ is isomorphic to the quasi-state space of a $C^*-$algebra.
\endroster
\endproclaim

\example{Remarks}
If $F_1$ and $F_2$ are closed faces of $C^*-$algebras, we say they are
{\it isomorphic} if there is a $0$-preserving affine isomorphism which is
also a (weak${}^*$) homeomorphism.  An intrinsic characterization of
$pAp$ was observed in \cite{11} (a portion of the proof of 3.5 for which
no originality was claimed): $pAp$ is the set of continuous affine
functionals vanishing at $0$ on $F(p)$.  With help of \cite{15} one can
find intrinsic characterizations of $QC(p)$ and $SQC(p)$.  One of the
consequences of
[7, 4.4, 4.5] is that $pA^{**}p$ is the bidual of
the Banach space $pAp$.  In \cite{14} we will give an intrinsic characterization
of $pM(A)p$.  Thus many questions concerning a closed face of a $C^*-$algebra
$A$ can be treated intrinsically, without knowing what $A$ is.

The $C^*-$algebra of 5 is determined only up to Jordan $*-$isomorphism.
\endexample

\demo{Proof}
$1\Rightarrow 2$:  Since $SQC(p)\subset pAp$ and $pA^{**}p$ is the bidual
of $pAp$, $SQC(p)$ is dense in $pAp$ in the weak Banach space topology.
Therefore $SQC(p)$ is norm dense in $pAp$.  But $SQC(p)$ is norm closed (since
it is a $C^*-$algebra, for example).

$2\Rightarrow 3\Rightarrow 4$: Obvious.

$4\Rightarrow 1$:  Let $a\in A_{\sa}$.  Then $papap\in pAp$.  Let
$(e_i)_{i\in D}$ be an approximate identity of $\text{her}(1-p)$. Then
$pa(1-e_i)ap\to papap$.  By Dini's theorem for continuous functions on
$F(p)$, this convergence is uniform.  Thus $\Vert pa(1-e_i-p)ap\Vert\to 0$,
$\Vert(1-e_i-p)^{1/2}ap\Vert\to 0$, and $\Vert(1-e_i-p)ap\Vert\to 0$.  It
follows that $(1-p)ap\in Ap$, since $Ap$ is closed by an argument similar to
[7, 4.4]. If $(1-p)ap=xp$ for $x$ in $A$, then $pxp=0$.  Therefore
$x\in L+R$, where $L=\{b\in A:bp=0\}$ and $R=L^*=\{b\in A:pb=0\}$, (proof of
[7, 4.4]). Since $Lp=0$, $(1-p)ap=rp$ for some $r$ in $R$.  Then if
$a'=a-r-r^*$, $pa'p=pap$ and $a'p=pa'$.  Thus $pap\in SQC(p)$.  This shows
2, but since $pAp$ is $\sigma$-weakly dense in $pA^{**}p$, it is obvious
that 2 implies 1. 

That 3 implies 5 is obvious from previous remarks and is also essentially
included in the proof of [7, 4.5].

That 5 implies 2 is also obvious from previous remarks and the fact
(\cite{{4}, Theorem III.3}) that 2 is true when $p=1$.
\enddemo

\proclaim{Theorem 3.2}
Let $A$ be a $\sigma$-unital $C^*-$algebra, $p$ a closed projection in
$A^{**}$, and let $B=SQC(p)$.  If $B$ is non-degenerately embedded in
$pA^{**}p$, then $M(B)$ is naturally isomorphic to $QC(p)$.
\endproclaim

\example{Remarks}
1.  When $B=pAp$ (i.e., when the conditions of 3.1 hold), this result was partly
proved in [7, 4.5].

2.  It follows from 3.2 that if $SQC(p)$ is non-degenerate in $pA^{**}p$ and
if $p$ does not satisfy (MSQC), then $p$ does not satisfy (MQC).  This is
so because $M(B)\subset B''$.
\endexample

\demo{Proof}
Let $A^{**}$ be represented on $H$ via the universal representation of $A$.
The non-degeneracy hypothesis means that $B$ is non-degenerately represented
on $pH$.  Therefore $M(B)$ is isomorphic to the idealizer of $B$ in
$B(pH)$.  It follows that if $F$ is a closed subset of ${\Bbb R}$ and $h$
is in $M(B)_{\sa}$ then there is a hereditary $C^*-$subalgebra $B_0$ of
$B$ such that any approximate identity of $B_0$ converges to $p-E_F(h)$,
where the spectral projection is computed in  $B(pH)$.  Let
$\overline B=\{a\in A: ap=pa\}$, and let $\overline B_0$ be the inverse image
of $B_0$ in $\overline B$.  If $q$ is the limit in $B(H)$ of an approximate
identity of $\overline B_0$, then $q$ is an open projection in $A^{**}$,
$qp=pq$, and $qp=p-E_F(h)$. Thus $E_F(h)$ is $p\wedge(1-q)$, a closed
projection in $A^{**}$, and $h$ is in $QC(p)$.

Conversely, if $x\in QC(p)$ and $b\in B$, then $x=p\overline x$ and
$b=p\overline b$ where $\overline x\in M(A)$, $\overline x p=p\overline x$, and
$\overline b\in\overline B$.  Then $xb=p\overline x\overline b\in B$ and $bx=p\overline b
\overline x\in B$. Thus $x\in M(B)$.
\enddemo

\demo{Remark}
The $\sigma$-unitality was used only in the second part of the proof.
\enddemo

\proclaim{Theorem 3.3}
If $A$ in 3.1 is $\sigma$-unital, then the following conditions are
equivalent to 1-5 of 3.1:
\roster
\item"6."
$pAp\subset QC(p)$.

\item"7."
$pM(A)p=QC(p)$.
\endroster
\endproclaim

\demo{Proof}
That $2\Rightarrow 6$ is obvious.

$6\Rightarrow 3$:  Let $x$ be in $pAp$ and let $a$ be in $A$.  Write
$x=p\overline x$ where $\overline x\in M(A)$ and $\overline x p=p\overline x$.  Then
$x pap=p\overline x pap=p^2\overline xap\in pAp$.

That 7 implies 6 is obvious.

$2\Rightarrow 7$:  Clearly we have the non-degeneracy required for 3.2.
Let $x$ be in $pAp$ and let $y$ be in $M(A)$.  Write $x=p\overline x$ where
$\overline x$ is in $A$ and $p\overline x=\overline x p$.  Then
$xpyp=p(\overline x y)p\in pAp$, and $pypx=p(y\overline x)p\in pAp$.  Thus, in the
notation of 3.2, $pyp\in M(B)$, and hence $pyp\in QC(p)$.
\enddemo

\example{Example 3.4}
In this example $p$ is closed, infinite rank, abelian, and atomic, and
$pA^*p$ is norm separable.  Also $SQC(p)=\{0\}$ but $p$ satisfies
$(MQC)$.  In particular, $p$ is a counterexample for the question raised
in Section 0 about isolated points.  In fact, if $p_0$ is a minimal projection,
$p_0\le p$, and $p-p_0$ is closed, then obviously $p_0\in SQC(p)$.

Let $A=C([0,1])\otimes\Cal K$.  Here $\Cal K$ is the algebra of compact
operators on a separable infinite dimensional Hilbert space $H$, 
$\{e_1,e_2,\dots\}$ is an orthonormal basis of $H$, and $P_n$ is the 
projection on $\text{span}\{e_1,\dots,e_n\}$.  A criterion for weak
semicontinuity from [11, \S5.G] will be used to describe closed projections
in $A^{**}$.  A closed projection is given by a projection-valued function
$P:[0,1]\to B(H)$ such that if $h$ is any weak cluster point of $P(y)$ as
$y\to x$, then $h\le P(x)$.  More precisely, $P$ describes the atomic part of a closed
projection $p$, and $P$ determines $p$ since a closed projection is
determined by its atomic part.  (In our case $p$ will equal its atomic
part.)  We will construct a countable subset $S$ of $[0,1]$ and unit
vectors $v(x)$ for each $x$ in $S$.  For $x$ in $S$, $P(x)$ is the rank one
projection on ${\Bbb C}v(x)$, and for $x$ not in $S$, $P(x)=0$.
\endexample

The following trivial lemma is stated for purposes of reference:

\proclaim{3.4.1}
Let $\{x_i\}$ be a sequence of distinct points in $[0,1]$ and let $D$ be a
countable subset of $[0,1]$.  Then there are distinct points $y_{ij}$ in
$[0,1]\setminus (\{x_i\}\cup D)$ such that $|y_{ij}-x_i|\le 2^{-(i+j)}$.
\endproclaim

We will take $S=\bigcup_0^{\infty}\ S_n$, a disjoint
union, where $S_n$ and $v|_{S_n}$ will be constructed recursively so that
$\Vert P_n v(x)\Vert\le n^{-\frac12}$ for $x$ in $S_n$.

\roster
\item""
Step 0:\quad Take $S_0=\{\dfrac12\}$, $v(\dfrac12)=e_1$.

\item""
Step 1:\quad Take $S_1=\{x_i\}$ where the $x_i$'s are distinct,
$x_i\not=\dfrac12$, and $x_i\to\dfrac12$ as $i\to\infty$.  Let
$v(x_i)=2^{-\frac12}e_1+2^{-\frac12}e_{i+1}$ for $i=1,2,\dots.$

\item""
$\vdots$

\item""
Step n ($n>1$, step $n-1$ already completed): Write $S_{n-1}=
\{x_1, x_2,\dots\}$. Choose $y_{ij}$'s as in 3.4.1 with $D=\cup_0^{n-2} S_k$.  
Let $S_n=\{y_{ij}:i,j=1,2,\dots\}$ and $v(y_{ij})=n^{-\frac12}v(x_i)+
(1-n^{-1})^{\frac12} w_{ij}$, where $w_{ij}$ is a unit vector such
that $(w_{ij}, v(x_i))=0$ and $P_{i+j+n}w_{ij}=0$.
\endroster

The first step in the proof is to show that we get a closed projection.
Thus we may assume given a sequence $(t_r)$ in $[0,1]$ such that
$t_r\to t$ and $P(t_r)\overset w\to\longrightarrow\ h$. We must show
$h\le P(t)$.  We have no difficulty if $P(t_r)=0$.  Thus we may assume, after
passing to a subsequence, that $t_r\in S_{n(r)}$.  If $n(r)\to\infty$, then
since $\Vert P_{n(r)}P(t_r)P_{n(r)}\Vert\le n(r)^{-1}$, we must have $h=0$.
Thus, after again passing to a subsequence, we may assume $n(r)=n$, $\forall r$.
Now it is easy to see by induction that $\bigcup_0^n\ S_k$ is closed.  In
fact, every cluster point of $S_n$ is in 
$\overline{S_{n-1}}=\bigcup_0^{n-1}\ S_k$.
The proof that $h\le P(t)$ will be left to the reader in the cases $n=0$, $n=1$.
If $n>1$, write $t_r=y_{i(r)j(r)}$ in the notation of step $n$.  If
$i(r)+j(r)\nrightarrow\infty$, we may assume, after passing to a subsequence, that
$t_r=t$, $\forall r$, a trivial case.   If $i(r)+j(r)\to\infty$, then
$t\in S_m$ for some $m<n$.  We use induction on $n-m$.  First suppose
$i(r)\nrightarrow\infty$.  Then, passing to a subsequence, we assume $i(r)=i$,
$\forall r$.  Then $t=x_i$, and the construction shows that $h=n^{-1}P(x_i)$.
If $i(r)\to\infty$, let $t_r'=x_{i(r)}$.  Then $t_r'\to t$, and
$[v(t_r)-n^{-\frac12}v(t_r')]\overset w\to\longrightarrow 0$.  This shows
that $h=n^{-1}h'$, where $P(t_r')\overset w\to\longrightarrow h'$.
Since $h'\le P(t)$ by induction, we conclude that $h\le P(t)$, as desired.

Now since $A$ is separable, every state in $F(p)$ is the resultant of a
probability measure on $F(p)\cap P(A)$.  Since $F(p)\cap P(A)$ is countable,
the integral is a Bochner integral and thus the resultant is an atomic
state.  This shows that $p$ is atomic, as claimed.  Also, $pA^*p$ is norm
separable, being isometrically isomorphic to $\ell^1(S)$.  That $p$ is
abelian, in other words $pA^{**}p$ is abelian, is now obvious (cf \cite{10}).

Now if $h$ is in $pA^{**}p$, $h$ is determined by a function $\lambda$ in
$\ell^\infty(S)$, where $h(x)=\lambda(x)P(x)$, $x\in S$.  If $h$ is in
$SQC(p)$, then $h=p\overline h$, where $\overline h\in A$ and $p\overline h=\overline h p$. In particular, $\overline h(\cdot)$ is a norm continuous function from
$[0,1]$ to $\Cal K$.  If $x\in S$, there is a sequence $(x_n)$ in $S$ such
that $x_n\to x$ and $P(x_n)\overset w\to\longrightarrow tP(x)$, where
$0<t<1$.  Since $P(\cdot)\overline h(\cdot)P(\cdot)=\lambda(\cdot)P(\cdot)$, we
conclude that $\lambda(x_n)\to t\lambda(x)$. Since also
$h^2\in SQC(p)$, we also have $\lambda(x_n)^2\to t\lambda(x)^2$.  This implies
$\lambda(x)=0$.  Since $x$ is arbitrary, $h=0$.  (The only property of $h$
actually used
is that $h$, $h^2\in pAp$.)

Finally we note that any continuous function on $[0,1]$ gives rise to an
element $\overline h$ of the center of $M(A)$.  Thus $p\overline h\in QC(p)$.  It is
easy to see that such elements of $QC(p)$ generate $pA^{**}p$ as a
$W^*-$algebra, and hence $p$ satisfies $(MQC)$.

\example{Example 3.5}
By modifying the previous example, we can obtain either of the following:
\roster
\item"(a)"
a compact projection $\tilde p$ such that $QC(\tilde p)={\Bbb C}\tilde p$

\item"(b)"
a closed projection $p_1$ such that $SQC(p_1)=\{0\}$ and $QC(p_1)={\Bbb C}p_1$.
\endroster

In both cases we will still have $p$ infinite rank, abelian, and atomic and
$pA^*p$ norm separable, and of course $A$ will still be separable.

(a)  Let $A$ and $p$ be as in 3.4, and consider $\widetilde A$ and
$\tilde p=p+p_\infty$.  $\widetilde A^{**}$ is identified with
$A^{**}\oplus{\Bbb C}$ and $p_\infty$ has its usual meaning, so that
$p_\infty=0\oplus 1$ in $A^{**}\oplus{\Bbb C}$ and $\tilde p=p\oplus 1$.
Then $\tilde p$ is closed, and hence compact, in $\widetilde A^{**}$.  Suppose
$x=\lambda 1_{\widetilde A}+a$, $\lambda\in{\Bbb C}$, $a\in A$, and $x\tilde p=\tilde p x$.
Then $ap=pa$, and hence by 3.4, $ap=0$.  It follows easily that
$\tilde px=\lambda\tilde p$.  Therefore $QC(\tilde p)={\Bbb C}\tilde p$.

(b)  We will use a $C^*-$algebra $A_1$ which is a maximal hereditary
$C^*-$subalgebra of $\widetilde A$.  Let $p_0$ be the minimal projection in
$\widetilde A^{**}$ (actually in $A^{**}$) corresponding to the projection
$P(\frac12)$ in the notation of 3.4 ($p_0$ corresponds to the pure state
$\varphi_0$ where $\varphi_0(a)=(a(\frac12)e_1,e_1)$).  Then $p_0\le p\le\tilde
p$.  Let $A_1=\text{her}(1-p_0)$ and $p_1=\tilde p-p_0$ in $A_1^{**}$.
Then $A_1^{**}$ is identified with $(1-p_0)\widetilde A^{**}(1-p_0)$.  It is easy
to see that $p_1$ is closed in $A_1^{**}$:  The complementary projection
to $p_1$ in  $A_1^{**}$ is $1-p_0-p_1=1-\tilde p$, and this supports a
hereditary $C^*-$subalgebra of $\widetilde A$ which happens to be contained in
$A_1$ also.  If $x$ is in $A_1$ and $xp_1=p_1x$, then $x$ is also in $\widetilde A$
and $x\tilde p=\tilde px$.  Thus by (a), $\tilde px=\lambda\tilde p$ and
hence $p_1x=\lambda p_1$.  But $x\in A_1$ implies $xp_0=p_0x=0$.  Since
$\tilde p x=\lambda\tilde p$ implies $p_0x=\lambda p_0$,  $\lambda=0$.
Therefore $SQC(p_1)=\{0\}$.

Now $A_1$ can be regarded as the set of all norm continuous functions
$f:[0,1]\to\widetilde{\Cal K}$ such that $f(\frac12)P(\frac12)=P(\frac12)f(\frac12)=0$
and the image of $f$ in $\widetilde{\Cal K}/\Cal K$ is constant.  Since $\frac12$
is not an isolated point of $[0,1]$, $M(A_1)$ can be regarded as a set of
functions $g:[0,1]\setminus\{\frac12\}\to\widetilde {\Cal K}$ (cf\. \cite{{7}, Theorem 3.3} and note
that $\widetilde {\Cal K}$ is unital).  The requirements on $g$ are:
\roster
\item"(i)"
$g$ is norm continuous and bounded.

\item"(ii)"
$\underset{t\to\frac12}\to\lim\ (1_{\widetilde{\Cal K}}-P(\frac12))g(t)(1_{\widetilde{\Cal K}}-P(\frac12))$ exists in
norm.

\item"(iii)"
$\underset{t\to\frac12}\to\lim\Vert P(\frac12)g(t)(1_{\widetilde{\Cal K}}-P(\frac12))\Vert=
\underset{t\to\frac12}\to\lim\Vert(1_{\widetilde{\Cal K}}-P(\frac12))g(t)P(\frac12)\Vert=0$.

\item"(iv)"
If we write $g(t)=\lambda(t)1_{\widetilde{\Cal K}}+x(t)$, $\lambda(t)\in{\Bbb C}$, $x(t)\in\Cal K$,
then $\lambda(\cdot)$ is a constant.
\endroster
To see these, the main thing to note is that the constant function
$1_{\widetilde{\Cal K}}-P(\frac12)$ is in $A_1$.

Now assume $g$, as above, commutes with $p_1$.  Then $x(t)$ commutes with
$P(t)$ for all $t$ in $S\setminus\{\frac12\}$, in the notation of 3.4.
Just as in 3.4, this implies $P(t)x(t)=0$ for $t$ in $S\setminus\{\frac12\}$;
i.e., $p_1x=0$ and $p_1g=\lambda p_1$.  Thus $QC(p_1)={\Bbb C}p_1$.
\endexample

\example{Example 3.6}
Here we show, by a simpler example, how badly Theorem 3.2 can fail when the
non-degeneracy hypothesis is eliminated.  By \cite{{7}, Theorem 2.7}, if $B$ is a non-unital
separable $C^*-$algebra, then $M(B)$ is non-separable.  Thus if $A$ is
separable and $SQC(p)$ is non-unital (in particular non-trivial), and if
the conclusion of 3.2 is true, then $QC(p)$ is much larger
than $SQC(p)$.  In this example, $SQC(p)$ is (infinite dimensional and)
non-unital and $QC(p)=SQC(p)+{\Bbb C}p$.

Let $A=c\otimes\Cal K$.  Thus $A^{**}$ can be identified with the set
of bounded collections $\{h_n:1\le n\le\infty\}$, $h_n\in B(H)$.  Let
$v_n=2^{-\frac12}e_1+2^{-\frac12}e_{n+1}$, $n<\infty$, $v_\infty=e_1$, let
$p_n$ be the projection with range ${\Bbb C}v_n$, and let $p=\{p_n\}$ in
$A^{**}$.  Then $p$ is closed since $p_n\overset w\to\longrightarrow
\frac12 p_\infty$, and clearly $p$ is abelian.  Any element of $pA^{**}p$ is
given by $h_n=\lambda_np_n$, $1\le n\le\infty$, $\{\lambda_n\}$ bounded.
An easy argument, which is part of 3.4, shows that if $h\in SQC(p)$ then
$\lambda_\infty=0$ and $\lambda_n\to 0$ as $n\to\infty$.  Conversely, any
such $h$ is in $SQC(p)$;  in fact $h\in A\cap pA^{**}p$. Thus
$SQC(p)\cong c_0$.  Next we show that $h\in QC(p)$ implies $\lambda_n\to\lambda_\infty$.  If this is false for $h=h^*$, then there is a closed subset $F$
of ${\Bbb R}$ such that $\lambda_\infty\notin F$ and $\lambda_n\in F$
for infinitely many $n$.  If $q=E_F(h)$, then $q_\infty=0$ and
$q_n=p_n$ for infinitely many $n$.  Since $p_n\to\frac12 p_\infty\not= 0$, 
$q$ is not closed and $h$ is not $q$-continuous on $p$.  Thus $QC(p)\cong c$
and $QC(p)/SQC(p)$ is one dimensional.
\endexample
\medskip
\noindent
{\bf 4.  Closed faces with  (NCEB).}

If $\widehat A$ is the spectrum of $A$ and $p$ is a projection in $A^{**}$, we
will denote by $X$ the set of all $[\pi]$ in $\widehat A$ such that
$\pi^{**}(p)\not= 0$.  For $[\pi]$ in $X$ let $p_{[\pi]}$ be the atomic
projection in $A^{**}$ corresponding to $\pi^{**}(p)$.  Thus
$z_{\at}p=\sum_{x\in X} p_x$.  If  $p$ is closed, or even universally
measurable, then $p$ is determined by the $p_x$'s.  If $\varphi$ and $\psi$
are in $(0,\infty)P(A)$, we will say that $\varphi$ and $\psi$ are
{\it equivalent}, and write $\varphi\sim\psi$, if the pure states
$\dfrac \varphi{\Vert\varphi\Vert}$ and $\dfrac\psi{\Vert\psi\Vert}$ are
equivalent.

The proof of the next theorem and some of the other geometric arguments
in this paper were inspired by Glimm \cite{16}.

\proclaim{Theorem 4.1}
If $p$ is a projection in $A^{**}$ and if $p$ satisfies (NCEB), then
$p_x$ is finite rank, $\forall x\in  X$.
\endproclaim

\demo{Proof}
Let $\pi$ be an irreducible representation belonging to $x$, and let
$H_x$ be the range of $\pi^{**}(p)$.  If the conclusion is false, there is
an infinite orthonormal sequence, $\{e_1,e_2,\dots\}$, in $H_x$.  Choose
$t>0$ such that $[P(A)\cap F(p)]^-\subset [t,1]P(A)\cup\{0\}$
and choose $s$ such that $0< s< t$.  Let 
$v_n=s^{1/2}e_i+(1-s)^{1/2}e_n$,
$n>2$, where $i$ is $1$ or $2$.  Define $\varphi_n$, $\psi_n$ in
$P(A)\cap F(p)$ by $\varphi_n(a)=(\pi(a)v_n,v_n)$, $\psi_n(a)=(\pi(a)e_n,e_n)$. 
Let $\theta$ be any cluster point of $(\psi_n)$ in $Q(A)$.  Since
$e_n\overset w\to\longrightarrow 0$, $(\pi(a)e_i,e_n)\to 0$, $\forall a\in A$.
Therefore $s\psi_i+(1-s)\theta$ is a cluster point of $(\varphi_n)$.  If
$\theta=0$, we have a contradiction to (NCEB), since $0<s<t$.  Therefore
$\theta\in[t,1]P(A)$.  Since we must also have $s\psi_i+(1-s)\theta\in
[t,1]P(A)$, it follows that $\theta=r_i\psi_i$ for some $r_i\ge t>0$. We have
shown that   $\theta=r_1\psi_1$ and $\theta=r_2\psi_2$, a contradiction.
\enddemo

For the rest of this section we assume that $p$ is closed and satisfies (NCEB).  Let\newline
$\widetilde X=[P(A)\cap F(p)]^-\setminus\{0\}$.  Then $\widetilde X\subset
F(p)\cap[t,1]P(A)$ and $\widetilde X$ is locally compact, since $\widetilde X\cup\{0\}$
is closed.  
We identify $X$ with the set of equivalence classes in
$\widetilde X$ via $f:\widetilde X\to X$, where $f(\varphi)=[\pi_\varphi]$.  Give
$X$ the quotient topology arising from $f$.

\proclaim{Lemma 4.2}
$f$ is a closed map.
\endproclaim

\demo{Proof}
The main point is to show the following:  If $(\varphi_i)_{i\in D}$ and
$(\psi_i)_{i\in D}$ are nets in $\widetilde X$ such that
$\varphi_i\sim\psi_i$, $\varphi_i\to\varphi$, and $\psi_i\to\psi$, then
either $\varphi=\psi=0$ or $\varphi,\psi\in\widetilde X$ and $\varphi\sim\psi$.
Assume this is false and consider first the case $\varphi=0$, $\psi\in\widetilde X$.
Let $\pi$ be the reduced atomic representation of $A$, $H=H_\pi$, and
choose vectors $u_i$, $v_i$ in $\pi^{**}(p)H$ of norm at most $1$ such that
$\varphi_i=(\pi(\cdot)u_i,u_i)$, $\psi_i=(\pi(\cdot)v_i,v_i)$.  If
$g_i(a)=(\pi(a)u_i,v_i)$, then $|g_i(a)|\le\Vert\pi(a)u_i\Vert=\varphi_i(
a^*a)^{1/2}\to 0$.  Therefore $g_i\to 0$.  Now choose $r_i$ in ${\Bbb R}$ such
that $\Vert w_i\Vert=1$, where $w_i=r_iu_i+(\frac t2)^{1/2}v_i$.
Since $\Vert u_i\Vert^2\ge t$, $\{r_i\}$ is bounded.  Let
$\theta_i=(\pi(\cdot)w_i,w_i)$.  It follows from the above that
$\theta_i\in F(p)\cap P(A)$ and $\theta_i\to\frac t2 \psi$.  Since
$0<\Vert\frac t2\psi\Vert <t$, this contradicts (NCEB).

Next assume $\varphi, \psi\in\widetilde X$ and $\varphi\not\sim\psi$.  Then there
are invariant subspaces $H_1$ and $H_2$ of $H$, corresponding to inequivalent
irreducible representations, and non-zero vectors $u$ in $H_1$, $v$ in $H_2$
such that $\varphi=(\pi(\cdot)u,u)$ and $\psi=(\pi(\cdot)v,v)$.  Let
$u_i$, $v_i$, and $g_i$ be as above with the extra condition that
$Re(u_i,v_i)\ge 0$.  Passing to a
subnet, we may assume $g_i\to g$, $g\in A^*$.  Since $|g_i(a)|\le
\varphi_i(a^*a)^{1/2}$, $\forall a\in A$, then $|g(a)|\le\varphi(a^*a)^{1/2}$.
From the Hahn-Banach and Riesz-Fisher theorems we see that
$g=(\pi(\cdot)u,u')$ for some $u'$ in $H$.  Clearly, we may assume
$u'\in H_1$.  Similarly, $|g_i(a)|\le\psi_i(aa^*)^{1/2}$, and hence
$|g(a)|\le\psi(aa^*)^{1/2}$.  Therefore $g=(\pi(\cdot)v',v)$ for some $v'$
in $H_2$. It follows that $g=0$ (\cite{20}).  Now choose $r_i$ in
${\Bbb R}_+$ such that $\Vert w_i\Vert=1$, where $w_i=r_i(u_i+v_i)$.  Since
$2t\le\Vert u_i+v_i\Vert^2\le 4$, $\{r_i\}$ is bounded and bounded away from
$0$.  If $\theta_i=(\pi(\cdot)w_i,w_i)$, then $\theta_i\in F(p)\cap P(A)$ and
every cluster point of $(\theta_i)$ is of the form $r^2(\varphi+\psi)$ for
some cluster point $r$ of $(r_i)$.  Since this last functional is not a
multiple of a pure state, this contradicts (NCEB).

To complete the proof of the lemma, we have to show that the saturation of a
closed set is closed.  Suppose $Y$ is a closed subset of $\widetilde X$,
$\varphi_i\in f^{-1}(f(Y))$, and $\varphi_i\to\varphi$ in $\widetilde X$.
Choose $\psi_i$ in  $Y$ such that $\varphi_i\sim\psi_i$.  Passing to a subnet
if necessary, we may assume $\psi_i\to\psi$.  By what has already been
proved $\psi\in\widetilde X$ and $\psi\sim\varphi$.  Since $Y$ is closed,
$\psi\in Y$ and hence $\varphi\in f^{-1}(f(Y))$.  Thus $f^{-1}(f(Y))$ is
closed (relative to $\widetilde X$).
\enddemo

\proclaim{Theorem 4.3}
$X$ is a locally compact Hausdorff space and $f$ is a proper map from
$\widetilde X$ to $X$.
\endproclaim

\demo{Proof}
The fibers of $f$, i.e., the sets $f^{-1}(\{x\})$, $x\in X$, are compact
(even norm compact) by 4.1.  This, 4.2, and the fact that $\widetilde X$
is locally compact Hausdorff imply that $X$ is locally compact Hausdorff, by
standard point set topology.  Any closed map with compact fibers is proper;
i.e., the inverse image of a compact set is compact.
\enddemo

\example{Remarks}
1.  It follows from 4.3, or it could be deduced directly from the proof of
4.2, that the saturation of a compact subset of $\widetilde X$ is compact.

2.  The topology of $X$ is stronger than, and in general unequal to, the
relative topology that $X$ inherits from the usual hull-kernel topology
of $\widehat A$.  In fact, using \cite{5} and 1.4, we can easily construct a
closed projection satisfying (NCEB) and even (CEB) such that $X$ is a countably
infinite discrete space and the image of $X$ in prim $A$ consists of one
point.  Thus the relative topology is trivial on $X$.
\endexample

\proclaim{Lemma 4.4}
Assume $p$ is an atomic closed projection satisfying (NCEB) and that $pA^*p$ is
norm separable.  Then for every closed subset $S$ of $X$, $\sum_{x\in S} p_x$
is a closed projection.
\endproclaim

\demo{Proof}
Since $pA^*p$ has a linear subspace isometric to $\ell^1(X)$, $X$ must be
countable.  Let $p_S=\sum_{x\in S} p_x$.  Then every element of $F(p_S)$
is the resultant of a probability measure supported by
$[F(p_S)\cap P(A)]\cup\{0\}$, and {\it a fortiori} supported by
$f^{-1}(S)\cup\{0\}$.  Since $f^{-1}(S)\cup\{0\}$ is compact, every element
of $F(p_S)^-$ is the resultant of a probability measure on $f^{-1}(S)\cup\{0\}$.
Since $f^{-1}(S)$ is the disjoint union of countably many fibers of $f$,
since each of these fibers is contained in $F(p_x)$ for some $x$ in $S$, and 
since each $p_x$ is finite rank and hence closed, it is easy to see that
any such resultant is in $F(p_S)$.  Thus $F(p_S)$ is closed and hence $p_S$ is
closed.
\enddemo

\proclaim{Corollary 4.5}
Under the same assumptions, if $p\not= 0$, there is a minimal projection
$p_0$ such that $p_0\le p$ and $p-p_0$ is closed.  Also for every non-zero
closed subprojection $p'$ of $p$, there is a minimal projection $p_0$ such
that $p_0\le p'$ and $p'-p_0$ is closed.
\endproclaim

\demo{Proof}
Since $X$ is countable and locally compact Hausdorff, the Baire category
theorem implies that $X$ has an isolated point $x_0$.  Let  $p_0$ be any
minimal subprojection of $p_{x_0}$.  Then $p-p_0=(p_{x_0}-p_0)+p_{X\setminus
\{x_0\}}$, the sum of two orthogonal closed projections.  Therefore $p-p_0$
is closed (\cite{{1}, Theorem II.7}).
\enddemo

\example{Remarks}
1.  In Section 6 we will generalize 4.4 and 4.5 by dropping the requirement
that $pA^*p$ be norm separable, but we will add the assumption that $A$ is
separable.  We are not sure what technical assumptions are really needed.

2.  Corollary 4.5 and Examples 3.4 and 3.5(a) constitute our results on the
``isolated point'' question raised in Section 0.  The second sentence of
4.5 is closely analogous to the definition of a scattered topological space
and less closely analogous to the definition of scattered $C^*-$algebras.
Obviously we have not found a necessary and sufficient condition for this
to hold.  Example 3.4 shows that we cannot replace (NCEB) by the weaker condition
$[F(p)\cap P(A)]^-\subset [0,1]P(A)$, and 3.5(a) shows we cannot weaken (NCEB)
to $[F(p)\cap P(A)]^-\subset\{0\}\cup [t,1]S(A)$.  Any closed face of a
a scattered $C^*-$algebra satisfies the conclusion of 4.5 but not necessarily
the hypothesis. Example 5.12 below, whose primary purpose is something
else, is a closed face satisfying the conclusion of 4.5 (the proof of this
is in 7.9), but not (NCEB), and which is not isomorphic to a closed face of
any scattered $C^*-$algebra.
\endexample

We now consider the geometry of $F(p)$ in more detail.

\proclaim{Theorem 4.6}
If $p$ is a closed projection satisfying (NCEB) and if $(x_i)_{i\in D}$ is a net
in $X$ converging to $x$, then there is a subnet $(x_j)_{j\in D}$ such that
one of the following holds: 

1.  We have rank $p_{x_j}=k\le n=\text{rank } p_x$,
$\forall j$;  and there are orthonormal bases $\{e_1^j,\dots,e_k^j\}$ of
range $\pi_j^{**}(p_{x_j})$ and $\{e_1,\dots,e_n\}$ of range
$\pi^{**}(p_x)$ and an $n\times k$ matrix $T$ such that $tI_k\le T^*T\le I_k$
and $\forall z\in{\Bbb C}^k$, $\varphi_j(z)\to\varphi(w)$, where $\pi_j$
and $\pi$ are irreducible representations belonging to $x_j$ and $x$,
$v_j=\sum_1^k\ z_me_m^j$, $v=\sum_1^n w_me_m$, $w=Tz$,
$\varphi_j(z)=(\pi_j(\cdot)v_j,v_j)$, and $\varphi(w)=(\pi(\cdot)v,v)$.

2.  There is $\varphi$ in $P(A)\cap F(p_x)$ such that every cluster point of
$(\varphi_j)$ is a multiple of $\varphi$, $\forall\varphi_j\in P(A)\cap 
F(p_{x_j})$.
\endproclaim

\demo{Proof}
If rank $p_{x_i}\nrightarrow\infty$, we first pick a subnet such that rank $p_{x_j}=k$,
$\forall j$.  If rank $p_{x_i}\to\infty$, we
must show there is a subnet satisfying 2;  and we do this by contradiction.
Thus assume there are a subnet $(x_j)$ and pure states $\theta_j$,
$\psi_j$ in $F(p_{x_j})$ such that $(\theta_j)$ and $(\psi_j)$ converge
to non-proportional elements of $F(p_x)$.  In the first case choose an
arbitrary orthonormal basis $\{e_1^j,\dots,e_k^j\}$ of
range $\pi^{**}(p_{x_j})$. In the second case let $k=n+1$ and choose an
orthonormal set $\{e_1^j,\dots,e_k^j\}$ in range $\pi^*(p_{x_j})$ such
that $\theta_j=(\pi_j(\cdot)v_j,v_j)$ and $\psi_j=(\pi_j(\cdot)v_j',v_j')$
with $v_j,v_j'$ unit vectors in span $\{e_1^j,\dots,e_k^j\}$.

In both cases define $f_{\ell m}^j$ in $A^*$ by $f_{\ell m}^j=(\pi_j(\cdot)e_m^j,e_\ell^j)$, $1\le \ell,m\le k$.  Passing to a subnet, we may assume
$f_{\ell m}^j\to f_{\ell m}$, $\forall \ell, m$.  Since the matrix
$[f_{\ell m}^j]$ represents a positive linear functional on $A\otimes M_k$,
the same must be true of the matrix $[f_{\ell m}]$.  The GNS representation of $A\otimes M_k$ induced by $[f_{\ell m}]$ must be of the form
$\tilde\pi\otimes id$ for some respresentation $\tilde\pi$ of $A$, and
$[f_{\ell m}]$ must be the vector state induced by a vector
$(u_1,\dots,u_k)$ in $H_{\tilde\pi}\oplus\dots\oplus H_{\tilde \pi}$. In
other words, $f_{\ell m}=(\tilde\pi(\cdot)u_m,u_\ell)$.  Since
$f_{\ell\ell}\in[t,1][P(A)\cap F(p_x)]$, $\tilde\pi\cong\pi\oplus\dots
\oplus\pi$. Thus we may write $u_\ell=(u_{\ell 1},\dots,u_{\ell r})$, $r\le k$, where $u_{\ell p}\in \text{ range}\ \pi^{**}(p_x)$.

Now $f_{\ell\ell}=\sum_1^r (\pi(\cdot)u_{\ell p},u_{\ell p})$. Since
$f_{\ell\ell}\in[t,1]P(A)$, there must be a non-zero vector $y_\ell$ in range
$\pi^{**}(p_x)$ such that $u_{\ell p}=\lambda_{\ell p}y_\ell$ with
$(\lambda_{\ell\cdot})$ a non-zero element of ${\Bbb C}^r$. If
$z\in{\Bbb C}^k$ and $\varphi_j(z)$ is as above, then $\varphi_j(z)=
\sum\bar z_\ell f_{\ell m}^j z_m$ and hence $\varphi_j(z)\to\sum\bar z_\ell f_{\ell m}z_m=(\tilde\pi(\cdot)\sum z_\ell u_\ell,\sum z_\ell u_\ell)$. Since
this functional is a multiple of a pure state, the vectors
$\sum z_\ell u_{\ell p}$, $1\le p\le r$, must be proportional. Suppose, for
example, that $y_1$ and $y_2$ are linearly independent.  Then the choice
$z=(1,1,0,\dots,0)$ shows that $(\lambda_{1\cdot})$ and $(\lambda_{2\cdot})$ are
proportional.  For $\ell>2$, $y_\ell$ cannot be a multiple of both
$y_1$ and $y_2$.  Therefore all $(\lambda_{\ell\cdot})$ are proportional.
 Changing notation, we may write $u_{\ell p}=\lambda_py_\ell$.
Then $\varphi_j(z)\to(\sum_1^r |\lambda_p|^2)(\pi(\cdot)\sum z_\ell y_\ell,
\sum z_\ell y_\ell)$.  Now choose any orthonormal basis of range $\pi^{**}(p_x)$
and let $T$ be the matrix of $z\to (\sum_1^r |\lambda_p|^2)^{1/2}\sum z_\ell 
y_\ell$.  Since $t\Vert z\Vert_2^2\le\Vert \lim\ \varphi_j(z)\Vert\le
\Vert z\Vert_2^2$, we must have $tI_k\le T^*T\le I_k$.  This implies $k\le n$
so that $1$ is proved.  The other alternative is that span$\{y_\ell\}$ is
one dimensional.  Then let $\varphi'=(\pi(\cdot)y_1,y_1)$ and $\varphi=
\frac{\varphi'}{\Vert\varphi'\Vert}$.  Since each $f_{\ell m}$ is proportional 
to $\varphi$, $(\varphi_j(z))$ converges to a multiple of $\varphi$,
$\forall z\in{\Bbb C}^k$, and more generally every cluster point of
$(\varphi_j(z_j))$ is a multiple of $\varphi$ for any bounded net $(z_j)$
in ${\Bbb C}^k$.  If $k=\rank\ p_{x_j}$, this proves 2.  In the original
second case, $\rank\ p_{x_j}\to\infty$, $k=n+1$, this proves the
contradiction that establishes 2.
\enddemo

We say that a $C^*-$algebra $A$ satisfies (CEB) or (NCEB) if the closed
projection 1 in $A^{**}$ satifies (CEB) or (NCEB).  In [17, \S5]
Glimm proved a necessary and sufficient condition for $A$ to satisfy a
property weaker than (NCEB), $\overline{P(A)}\subset[0,1]P(A)$.  His
condition is:
\roster
\item"(i)"
$A$ is CCR,
\item"(ii)"
$\widehat A$ is Hausdorff, and
\item"(iii)"
$[\pi]\in\widehat A$ and $\text{dim}\ \pi>1$ implies $[\pi]$ is regular.
\endroster
Given (i) and (ii), (iii) can be restated as follows:  If I is the ideal
of $A$ such that $\widehat I=\{[\pi]\in\widehat A: \text{dim}\ \pi>1\}$, then $I$ is a
continuous trace $C^*-$algebra.  (See [27] for the theory of continuous trace 
$C^*-$algebras.)  It is presumably an easy exercise to
derive a characterization of $C^*-$algebras satisfying (CEB) or (NCEB) (they
are equivalent for $C^*-$algebras) from Glimm's result.  In Corollary 4.7
below we derive such a characterization instead from 4.1-4.6.  The purpose
is not to put this result on the record, so long after \cite{17}.  The
purpose is as follows:  The class of closed faces of $C^*-$algebras admits
more varied behavior than the class of $C^*-$algebras. One illustration of this
is the contrast between the facts on the isolated point question for atomic
closed faces of $C^*-$algebras and the facts on scattered $C^*-$algebras
(\cite{18}, \cite{19}).  Another illustration is the contrast between 4.6
and 4.7. (We will show by example that all of the behavior contemplated by
4.6 really occurs.)  The exercise of deriving 4.7 from 4.1 to 4.6 gives some
insight into why the behavior of closed faces is more varied than that of
$C^*-$algebras.

If $A$ is a CCR $C^*-$algebra with Hausdorff spectrum, then $A$ is isomorphic
to the set of continuous sections vanishing at $\infty$ of a continuous
field, $\Cal A(x)$, $x\in\widehat A$, of elementary $C^*-$algebras.  If
$x_0\in \widehat A$ and $\Cal A(x_0)$ is one dimensional, then there is a
continuous section $e(\cdot)$ such that $e(x_0)=1_{\Cal A(x_0)}$ and $e(x)$ is a
projection for $x$ in some neighborhood of $x_0$ (\cite{17}).  We will say
$A$ is {\it locally unital} at $x_0$ if $e(x)=1_{\Cal A(x)}$ in some neighborhood of
$x_0$.

\proclaim{Corollary 4.7}
The following are equivalent for a $C^*-$algebra $A$:
\roster
\item"1."
$A$ satisfies (CEB)
\item"2."
$A$ satisfies (NCEB)
\item"3." 
{\rm (i)}\ \ Every irreducible representation of $A$ is finite dimensional,
\item""
{\rm (ii)}\ \ $\widehat A$ is Hausdorff,
\item""
{\rm (iii)}\ \ $\forall n>1$, $\{[\pi]:dim\ \pi=n\}$ is an open subset of $\widehat A$, and
\item""
{\rm (iv)}\ \ $A$ is locally unital at each $[\pi]$ with $dim\ \pi=1$.
\endroster
\endproclaim

\example{Remark}
Condition 3(iii) says that the ideal $I$ discussed above is the $c_0$ direct sum of
$n$-homogeneous $C^*-$algebras for various values of $n$. Thus the comparison
of 3 with Glimm's condition is clear.
\endexample

\demo{Proof}
$2\Rightarrow 3$:  (i) follows from 4.1 with $p=1$.  Since $p=1$, $X=\widehat A$.
Since the map from $P(A)$ to $\widehat A$ is open for the hull-kernel
topology (\cite{17}), the hull-kernel topology is the quotient topology; i.e.,
our topology on $X$ agrees with the usual one when $p=1$.
Thus (ii) follows from 4.3.  Again since the map from
$P(A)$ to $\widehat A$ is open, if $\dim\ \pi>1$ and $[\pi_i]\to[\pi]$, then
after passing to a subnet, we can find $\varphi_i$, $\psi_i$ in $P(A)$ 
associated to $\pi_i$ such that the nets $(\varphi_i)$, $(\psi_i)$ converge
to distinct pure states associated to $\pi$.  Thus alternative 2 of 4.6
cannot hold, and $\lim\sup(\dim\ \pi_i)\le \dim\ \pi$.  It is always true
in a $C^*-$algebra that $\lim\inf(\dim\ \pi_i)\ge \dim\ \pi$ (but for a closed
face we can have $\lim\inf(\rank\ p_{x_i})< \rank\ p_x)$.  This shows (iii).
If $x_0$, $e(\cdot)$ are as above and $A$ is not locally unital at $x_0$,
then we can find $(x_i)$ such that $x_i\to x_0$ and $e(x_i)\not= 1$,
$\forall i$.  Then we can find $\varphi_i$ in $P(A)$ associated to $x_i$ such
that $\varphi_i(e)=\frac t2$.  It follows that $\Vert\varphi\Vert=\frac t2$
for any cluster point $\varphi$ of $(\varphi_i)$, in contradiction to
(NCEB).  This proves (iv).

That 1 implies 2 is obvious, and the proof that 3 implies 1 is left to the
reader.
\enddemo

\example{Examples 4.8}
(a) We can illustrate alternative 1 of 4.6 with $A=c\otimes \Cal K$.
Choose $k$ and $n$ with $k\le n$, $t>0$, and an $n\times k$ matrix $T$ such
that $tI_k\le T^*T\le I_k$. Let $S=(1-T^*T)^{1/2}$, a $k\times k$ matrix.  Let
$p_\infty$ be the projection on $\text{span}\{e_1,\dots,e_n\}$ and for $j<\infty$
let $p_j$ be the range projection of $\left(\smallmatrix T\\  S\endsmallmatrix
\right)$, where the matrix is regarded as a linear isometry from ${\Bbb C}^k$
to $\text{span}\{e_1,\dots,e_n,e_{n+j},\dots,e_{n+j+k-1}\}$. If
$p=\{p_j:1\le j\le\infty\}$, then $p$ is a closed projection in $A^{**}$,
$p$ satisfies (NCEB) ((CEB) if $t=1$) and 4.6.1 holds with the given matrix
$T$.  (Here we think of $x_j$ as $j$ and $x$ as $\infty$, and 
$\{e_1^j,\dots,e_k^j\}$ corresponds to the columns of $\left(\smallmatrix
T\\ S\endsmallmatrix\right)$.)

If we want a more complicated example, say one where two different
subsequences give two different matrices, we can easily modify the above.
Choose $k'\le n$ and an $n\times k'$ matrix $T'$ such that
$tI_{k'}\le {T'}^*T'\le I_{k'}$.  Let $\tilde p_{2j-1}$ be the above $p_j$, and let
$\tilde p_{2j}$ be the above $p_j$ constructed from $T'$ instead of $T$.

(b) As a first example for alternative 2 of 4.6, consider
$A_1=\{(a_n)_1^\infty:a_n\in\widetilde{\Cal K}\text{ and } (a_n)$ 
$\text{converges to a
scalar in norm}\}$.  Then $A_1^{**}$ can be identified with the set of
bounded collections $\{h_n:1\le n\le\infty\}$  such that $h_n\in B(H)\oplus
{\Bbb C}$ for $n<\infty$ and $h_\infty\in{\Bbb C}$.  Choose any sequence
$(n_j)$ of positive integers and define a closed projection $p$ in
$A_1^{**}$ by: $p=\{p_j\}$, $p_\infty=1_{\widetilde {\Cal K}}$, and $p_j$ is a rank $n_j$ projection
in $B(H)$ for $j<\infty$.  It is easy to see that $p$ satisfies (CEB) and
4.6.2.  This easy example shows that there is no restriction on rank $p_{x_j}$
when 4.6.2 holds, but this is all that it shows.

(c)  For more complicated examples, in particular examples where some
subsequences satisfy 4.6.1 and others 4.6.2, we can use $A_2=A_1\otimes
\Cal K$.  Then $A_2^{**}\cong A_1^{**}\overline\otimes B(H)$.

The construction in (a) above can also be used for $A_2$.  Let
$\tilde p_\infty=1\otimes p_\infty$ and $\tilde p_j=q_0\otimes p_j$ for
$j<\infty$, where the $p_j$'s are as in (a) and $q_0$ is a rank 1 projection
in the $B(H)$-component of $\widetilde{\Cal K}^{**}$.  It is easy to see that $\tilde p$
is closed in $A_2^{**}$ and that $F(\tilde p)$ is isomorphic to the closed face
$F(p)$ of $(c\otimes\Cal K)^{**}$.

We can also construct examples of 4.6.2 using $A_2$.  Let $T$ be a
positive $k\times k$ matrix such that $tI_k\le T^2\le I_k$, and let $u$ be a unit
vector in $\text{span}\{e_1,\dots,e_n\}$ where $k$ and $n$ are arbitrary.  Let
$S=(1-T^2)^{1/2}$ and define a closed projection $\tilde p$ in $A_2^{**}$ by:
$\tilde p=\{\tilde p_j:1\le j\le\infty\}$, $\tilde p_\infty=1\otimes p_\infty$
for $p_\infty$ the projection on $\text{span}\{e_1,\dots,e_n\}$, and 
$\tilde p_j$ is 
the range projection of $\left(\smallmatrix T\\  S\endsmallmatrix\right)$
where now $\left(\smallmatrix T\\  S\endsmallmatrix\right)$ sends
${\Bbb C}^k$ to
$$
\text{span}\{e_1\otimes u, e_2\otimes u,\dots,e_k\otimes u, e_1\otimes e_{n+j},\dots,
e_1\otimes e_{n+j+k-1}\}.
$$
Then 4.6.2 holds with $\varphi$ given by $\varphi(a)=(a_\infty u,u)$.  Also
the columns of $\left(\smallmatrix T\\  S\endsmallmatrix\right)$ give an
orthonormal basis $\{e_1^j,\dots,e_k^j\}$ of range $\tilde p_j$, and,
using the notation of 4.6.1, $\varphi_j(z)\to\Vert Tz\Vert^2\varphi$.  It
is easy to see that $\tilde p$ satisfies (NCEB).

By using the idea of the second paragraph of (a), we can construct a closed
projection such that different subsequences exhibit different behavior.
Some subsequences can satisfy 4.6.1, with different choices of $T$ and $k$,
and some can satisfy 4.6.2 with $\varphi_j(z)\to\Vert Tz\Vert^2\varphi$ for
different choices of $T$, $k$, and $\varphi$.
\endexample

\example{Remark}
In 4.6.2 we showed only that every cluster point of $(\varphi_j)$ is a
multiple of $\varphi$ and did not describe which multiples arise.  When
rank $p_{x_j}$ is bounded, the same methods can easily be used to construct
a subnet and a positive $k\times k$ matrix $T$ such that $tI_k\le T^2\le I_k$
and $\varphi_j(z)\to\Vert Tz\Vert^2\varphi$.
\endexample
\medskip

\noindent
{\bf 5.  Type I Closed Faces and Atomic Closed Faces}

If $p$ is a projection in $A^{**}$, we say that $p$ or $F(p)$ is
{\it type I} if $pA^{**}p$ is a type I von Neumann algebra. Clearly $p$
is type I if and only if $c(p)$, the central cover of $p$, is type I.  Now
$F(p)$ is the normal quasi-state space of $pA^{**}p$, and for $\varphi$
in $F(p)$ the kernel of $\pi_\varphi$ contains $(1-c(p))A^{**}$.
Therefore $p$ is type I if and only if $\pi_\varphi$ is a type I representation
for all $\varphi$ in $F(p)$.  (It doesn't matter whether we look
at $\pi_\varphi$ or $\pi_\varphi^{**}$.)
Because $z_{\at}A^{**}$ is a type I $W^*-$algebra, every atomic
projection is type I.  
However, if we also require that $p$ be closed, or 
just universally measurable (say), it seems that the property of being type
I may be useful.

\proclaim{Lemma 5.1}
Let $A$ be a separable $C^*-$algebra and $p$ a type I closed projection in
$A^{**}$.  Let $\mu$ be a probability measure on $F(p)$, and let
$\pi=\int^\oplus \pi_\omega d\mu(\omega)$, the direct integal.  Then $\pi$
is a type I representation.
\endproclaim

\demo{Proof}
Let $\varphi=\int \omega d\mu(\omega)$, the resultant of $\mu$.  Then
$\varphi\in F(p)$, since $F(p)$ is closed.  Therefore $\pi_\varphi$ is type I,
and $\pi_\varphi$ is a subrepresentation of $\pi$. We claim $\pi$  and $\pi_\varphi$
have the same central support in $A^{**}$ (i.e. $\pi$ is quasi-equivalent
to $\pi_\varphi$).  Therefore $\pi$ is also type I.

To see the claimed quasi-equivalence, let $v_\omega$ be the cyclic vector
in $H_\omega$ produced by the GNS construction, and let $v=\int^\oplus
v_\omega d\mu(\omega)$, a vector in $H_\pi$.  Then $(\pi(a)v,v)=\varphi(a)$.
For every $\mu$-measurable subset $S$ of
$F(p)$ ($\mu$ is a Borel measure, and ``$\mu$-measurable''
means measurable with respect to the completion of $\mu$)
there is a projection $P_S$ in $\pi(A)'$ such that the corresponding
subrepresentation of $\pi$ is $\int_S^\oplus \pi_\omega d\mu(\omega)$.  It
is easy to see that $H_\pi$ is the smallest closed invariant subspace
containing $P_Sv$ for all such $S$.  Moreover the cyclic subrepresentation
of $\pi$ generated by $P_Sv$ is equivalent to a subrepresentation of 
$\pi_\varphi$.  These remarks complete the proof.
\enddemo 

The main fact needed from direct integral theory is
something that the author learned from G. W. Mackey and is expressed as a
lemma.  For the ideas in the proof see \cite{23},
pages 112-117, and \cite{24}, especially page 159.  The basic point
is that the direct integral decomposition into irreducibles of a type I
representation is almost unique.

\proclaim{Lemma 5.2 (Mackey)}
Let $A$ be a separable $C^*-$algebra, let $\pi_\cdot'$ and $\pi_\cdot''$
be measurable fields of irreducible representations of $A$ defined over
standard measure spaces $S'$ and $S''$, and let $\pi'=\int_{S'}^\oplus
\pi_s'd\mu'(s)$, $\pi''=\int_{S''}^\oplus \pi_s''d\mu''(s)$.   Assume that
$\pi_{s'}'$ is inequivalent to $\pi_{s''}''$, $\forall s'\in S'$,
$\forall s''\in S''$ and that $\pi'$ and $\pi''$ are  type I representations.
Then $\pi'$ and $\pi''$ are disjoint (i.e., their
central supports in $A^{**}$ are orthogonal).
\endproclaim

\proclaim{Lemma 5.3}
Let $A$ be a separable $C^*-$algebra and $p$ a type I closed projection in
$A^{**}$.  Assume $\mu$ and $\nu$ are positive finite measures on
$F(p)\cap P(A)$ such that $\int \omega d\mu(\omega)=
\int \omega d\nu(\omega)$.  Let $E$ be a saturated Borel
subset (or, more generally, a saturated $(\mu+\nu)$-measurable subset)
of $F(p)\cap P(A)$.  Then $\int_E \omega d\mu(\omega)=
\int_E \omega d\nu(\omega)$ and in particular $\mu(E)=\nu(E)$.
\endproclaim

\demo{Proof}
Let $\varphi=\int \omega d\mu(\omega)=\int \omega d\nu(\omega)$,
$\pi'=\int^\oplus \pi_\omega d\mu(\omega)$, and $\pi''=\int^\oplus
\pi_\omega d\nu(\omega)$.  As in the proof of 5.1, there are vectors $v'$
in $H_{\pi'}$ and $v''$ in $H_{\pi''}$ which induce the functional $\varphi$.
Thus there is a partial isometry $U$ which intertwines $\pi'$ and
$\pi''$ such that $v'$ is in the initial space of $U$ and $U{v'}=v''$.

Let $P_E'$ and $P_E''$ be the projections in $\pi'(A)'$ and $\pi''(A)'$
defined from $E$.  Thus $\mu(E)=(P_E'v',v')$ and $\nu(E)=(P_E''v'',v'')$.  By
5.2 and 5.1, $(1-P_E'')UP_E'=P_E''U(1-P_E')=0$.  Therefore $P_E'v'$ is in
the initial space of $U$ and $UP_E'v'=P_E''v''$.  The conclusion follows.
\enddemo

\proclaim{Lemma 5.4}
Let $A$ be a separable $C^*-$algebra, $p$ a type I closed projection in
$A^{**}$, and $E$ a saturated Borel subset of $F(p)\cap P(A)$. Then there is
a projection $p_E$ in $A^{**}$ such that $p_E\le p$ and $F(p_E)$ is the set
of resultants of probability measures on $E\cup\{0\}$.
\endproclaim

\demo{Proof}
Let $F_1$ be the set of resultants of probability measures on $E\cup\{0\}$.
We claim that $F_1$ is a norm closed sub-face of $F(p)$.  The result then
follows from \cite{{15}, Theorem 4.4 and p\. 396} (cf\. [25, 3.6.11]). 

To see the claim, note that by Choquet theory every element of $F(p)$ is the
resultant of a probability measure on $[F(p)\cap P(A)]\cup\{0\}$.  Let
$E'=[F(p)\cap P(a)]\setminus E$.  Then $F_1=\{\int \omega d\mu(\omega):\mu(E')=0\}$.  Suppose $\varphi_i=\int \omega d\mu_i(\omega)$, $i=1,2,$ and
$t\varphi_1+(1-t)\varphi_2\in F_1$, $0<t<1$.  By 5.3, $t\mu_1(E')+(1-t)\mu_2(E')=0$.  Therefore $\mu_1(E')=\mu_2(E')=0$, and $\varphi_1,\varphi_2\in F_1$.
Thus $F_1$ is a face.

To see that $F_1$ is norm closed, assume $\varphi=\int \omega d\mu(\omega)$
where $\mu(E')=\delta>0$.  We claim that $\dist(\varphi,F_1)\ge\delta$.
Suppose $\psi=\int\omega d\nu(\omega)$ where $\nu(E')=0$ and
$\Vert\varphi-\psi\Vert=r$.  Then $\varphi-\psi=\lambda_1-\lambda_2$
where $\lambda_1,\lambda_2\ge 0$ and $\Vert\lambda_1\Vert+
\Vert\lambda_2\Vert=r$.  If $\lambda_i=\int \omega d\mu_i(\omega)$, for
positive measures $\mu_1$, $\mu_2$ on $F(p)\cap P(A)$, then $\mu+\mu_2$
and $\nu+\mu_1$ have the same resultant.  Therefore by 5.3, 
$\mu(E')+\mu_2(E')=\nu(E')+\mu_1(E')$.  Therefore $\mu(E')\le\mu_1(E')\le r$.  
This proves the claim and completes the proof of the lemma.
\enddemo

\example{Remarks}
Although the conclusion of 5.4 has what we need, more is true.  Also
$F(p_E)\cap S(A)$ is a split face of $F(p)\cap S(A)$, the complement being
$F(p_{E'})\cap S(A)$.  This means that $p_E$ and $p_{E'}$ are centrally disjoint
projections and $p_E+p_{E'}=p$.  Also $p_E$ satisfies the barycenter
formula.  (The barycenter formula is discussed below before 5.13).  A
related statement is that $F(p_E)$ is closed under resultants.
The hypotheses of 5.4 could be weakened.  We could assume that $p$ satisfies
the barycenter formula instead of that $p$ is closed, and we could
assume $E$ universally measurable instead of Borel.
\endexample

\proclaim{Lemma 5.5}
Let $A$ be a separable $C^*-$algebra and $p$ a closed projection in
$A^{**}$.  If $\pi^{**}(p)$ has finite rank for every irreducible
representation $\pi$ of $A$, then $p$ is type I.
\endproclaim

\demo{Proof}
Let $\pi=\int^\oplus \pi_sd\mu(s)$ be a standard direct integral, where
each $\pi_s$ is irreducible.  Since $p$ is closed, $\pi^{**}(p)=
\int^\oplus \pi_s^{**}(p)d\mu(s)$, where $\pi_\cdot^{**}(p)$ is a Borel
operator field.  Therefore rank $(\pi_\cdot^{**}(p))$ is a Borel function,
and by hypothesis it is everywhere finite-valued.

From the above it follows that any representation of $A$ in a separable
Hilbert space can be written as a direct sum, $\pi=
\bigoplus_0^{\infty}\ \pi_n$, such that $\pi_n=\int_{S_n}^\oplus \pi_sd\mu(s)$ and
rank $(\pi_s^{**}(p))=n$, $\forall s\in S_n$.  It was shown by A. Amitsur
and J. Levitzki in [9] that there
is a non-commutative polynomial $G_n$ of $2n$ variables such that $G_n$ vanishes
on $M_n^{2n}$ but not on $M_{n+1}^{2n}$ (cf\. [21, Lemma 2], where a 
weaker but adequate result is proved).  Also if $G_n$ vanishes on $M^{2n}$
for a $W^*-$algebra $M$, then $M$ is a direct sum of type $I_k$ algebras
for $k\le n$.  Clearly $G_n$ vanishes on $[\pi_n^{**}(p)\pi_n(A)\pi_n^{**}(p)]^{2n}$, $n>0$, and hence by strong continuity $G_n$ vanishes on
$[\pi_n^{**}(pA^{**}p)]^{2n}$.  Therefore $\pi_n^{**}(pA^{**}p)$ is type I,
$\forall n$.  (For $n=0$, $\pi_0^{**}(p)=0)$.  If $z_n$ is the central support
of $\pi_n$ in $A^{**}$, and $z(\pi)$ is the central support of $\pi$, then
$z(\pi)=\sup_n\ z_n$.  Since we have shown that
$z_npA^{**}p$ is type I, $\forall n$, then $z(\pi)pA^{**}p$ is type I.
Since $\sup\{z(\pi):\pi\ \text{ as above}\} =1$, $pA^{**}p$ is type I.
\enddemo

\proclaim{Corollary 5.6}
If $A$ is a separable $C^*-$algebra, $p$ is a closed projection in $A^{**}$,
and if $p$ satisfies (NCEB), then $p$ is type I.
\endproclaim

\demo{Proof}
Combine 4.1 and 5.5.
\enddemo

We have already defined the concept of an atomic projection in $A^{**}$.  We
say that $p$ is {\it strongly atomic} if $p$ is atomic and $pA^*p$ is norm 
separable.  If $A$ is separable the separability of $pA^*p$ can be rephrased:
There are only countably many points $[\pi]$ in $\widehat A$ such that 
$\pi^{**}(p)\not=0$.

\example{Question 5.7}
If $A$ is separable, is every closed atomic projection in $A^{**}$ strongly
atomic?
\endexample

If $p$ is closed and atomic and if $\mu$ is a probability measure on
$F(p)\cap P(A)$, then $\int \omega d\mu(\omega)$ is in $F(p)$ and hence is an
atomic state.  If $A$ is separable, it follows from 5.3, for example,
that $\mu$ is supported by the union of countably many equivalence classes.
If $p$ is not strongly atomic, this means that there are uncountably many
equivalence classes in $F(p)\cap P(A)$ but every finite measure is concentrated
on the union of countably many.  It follows that the relation of equivalence of
pure states is not countably separated on $F(p)\cap P(A)$.  (If it were
countably separated, the quotient Borel space would be an uncountable
analytic Borel space (\cite{{24}, Theorem 5.1}) and hence would support a continuous
measure.  This measure could be lifted to $F(p)\cap P(A)$ by the von Neumann
selection lemma.)  In particular  $A$ is not type I.  Also $p$ does not
satisfy (NCEB), since the space $X$ of Section 4 is second countable and
hence countably separated when $A$ is separable.  This reasoning suggests
the following:

\example{Question 5.8}
If $A$ is a separable $C^*-$algebra and $p$ is a type I closed projection
in $A^{**}$, is equivalence of pure states countably separated on
$F(p)\cap P(A)$?
\endexample

Obviously 5.8 is analogous to Mackey's conjecture ([23, p. 85] or [24, p. 163]),
which was proved by Glimm in \cite{17}.  Of course \cite{17} proved much
more than Mackey's conjecture.  We do not know whether 
there is a structure theorem for type I closed faces of similar power to
Glimm's theorem.  Because the variety of closed faces of $C^*-$algebras is
so great, there is not enough evidence to support a conjecture on any of
these questions.

If the answer to 5.8 is yes for a particular $p$, then a standard form for 
elements of $F(p)\cap S(A)$ can be established.  Let $X$ be the set of
equivalence classes of $F(p)\cap P(A)$, an analytic Borel space which is in
one-one corespondence with a subset $\{[\pi_x]: x \in X\}$ of $\widehat A$.
Then an element $\varphi$ of $F(p)\cap S(A)$ is determined by a probability
measure $\mu$ on $X$ and a measurable function $f:X \to S(A)$ such that
$f(x)$ is supported by $\pi_x^{**}(p)$.
In fact $\varphi$ is the resultant of a probability measure $\overline \mu$
on $F(p)\cap P(A)$.  Even though $\overline \mu$ is not unique, 5.3 implies
its pushforward to $X$ is unique.  The function $f$ is obtained by writing
$\overline \mu =\int_X \nu_x d\mu(x)$, where $\nu_x$ is supported on the
equivalence class $x$, and $f(x)=\int \omega d\nu_x(\omega)$.  It can be
shown that $f$ is unique modulo null sets.  Thus, under the hypotheses
given, the Choquet decomposition of $\varphi$ is almost unique in a sense
roughly analagous to Mackey's result that the direct integral decomposition
of a type I representation into irreducibles is almost unique.

There is a converse question to 5.7 which we can answer.  The proof is valid
even for $A$ nonseparable.

\proclaim{Proposition 5.9}
If $A$ is a $C^*-$algebra and $p$ is a closed projection in $A^{**}$ such that
$z_{\at}pA^*p$ is norm separable, then $p$ is atomic and hence strongly
atomic.
\endproclaim

\demo{Proof}
There is an increasing sequence $(p_n)$ of finite rank projections such
that $p_n\to z_{\at}p$.  By 4.5.12 or 4.5.15 of \cite{25}, $z_{\at}p$ is
universally measurable.  Since $(1-z_{\at})p$ is a universally measurable
operator whose atomic part is $0$, $(1-z_{\at})p=0$ ([25, 4.3.15]).
\enddemo

The following lemma, or the ideas in its proof, might be useful in
connection with questions 5.7, 5.8.  It will also be used to prove a
complement to Glimm's theorem.

\proclaim{Lemma 5.10}
If $p$ is an atomic projection in $A^{**}$ such that $pA^*p$ is norm
separable, then $F(p)\cap P(A)$ is an $F_\sigma$ set relative to $P(A)$.
\endproclaim

\demo{Proof}
The lemma can be rephrased more concretely:  Let $\pi:A\to B(H)$ be an
irreducible representation, let $H_0$ be a separable closed subspace of
$H$, and let $P_0=\{(\pi(\cdot)v,v):v \text{ is a unit vector in $H_0$}\}$.
Then $P_0$ is an $F_\sigma$ set relative to $P(A)$.

The proof is similar to that of 4.1.  Let $H_1,H_2,\dots$ be an increasing
sequence of finite dimensional subspaces such that $H_0=(\cup_1^\infty H_n)^-$,
and let $p_n$ be the projection on $H_n$.  Let $V_n=\{v\in H_0:\Vert v\Vert=1
\text{ and } \Vert p_nv\Vert\ge \tfrac 12 \}$ and $P_n=\{(\pi(\cdot)v,v):
v\in V_n\}$.  Then $P_0=\cup_1^\infty P_n$, and we will show $P_n$ closed
relative to $P(A)$.  Suppose $v_i\in V_n$, $\varphi_i=(\pi(\cdot)v_i,v_i)$, and the net $(\varphi_i)$ converges to a pure state $\varphi$.  Passing to a
subnet if necessary, we may assume $v_i\overset w\to\longrightarrow v$
for some $v$ in $H_0$.  Clearly $\Vert v\Vert\le 1$ and 
$\Vert p_nv\Vert\ge \tfrac 12$.  Then $v_i=u_i+w_i$, where $u_i\to v$ in norm,
$w_i\overset w\to\longrightarrow 0$, and $(u_i,w_i)=0$.  Therefore
$(\pi(a)u_i,w_i)\to 0$, $\forall a\in A$.  Passing to a further subnet, we
may assume $(\pi(\cdot)w_i,w_i)$ converges to some $\psi$ in $Q(A)$.  Then
$\varphi=(\pi(\cdot)v,v)+\psi$.  Since $\varphi$ is pure, $\psi$ must be
proportional to $(\pi(\cdot)v,v)$.  Therefore $\varphi=(\pi(\cdot)v_1,v_1)$
where $v_1=v/\Vert v\Vert$.  Since $\Vert p_nv_1\Vert\ge\Vert p_nv\Vert$,
$\varphi\in P_n$.
\enddemo

\proclaim{Proposition 5.11}
If $A$ is a separable
$C^*-$algebra and $\pi:A\to B(H)$ is an irreducible representation such
that $\pi(A)\not\supset\Cal K(H)$, then there are uncountably many
inequivalent irreducible representations of $A$ with the same kernel as
$\pi$.
\endproclaim

\example{Remark}
Glimm's theorem
implies that there are uncountably many irreducibles with the same kernel,
but so far as we know, it was not previously known that that kernel can be
taken to be the same as the kernel of the
given $\pi$.
\endexample

\demo{Proof}
By replacing $A$ with its quotient by the kernel of $\pi$, we may reduce to
the case $\pi$ faithful.  Assume that $A$ has only countably many faithful
irreducible representations. Since $\widehat{A}$ is second countable, there
is a countable set $\{I_n\}$ of non-zero (closed, two-sided) ideals such that 
every non-faithful representation of $A$ vanishes on some $I_n$.  Then since
$[\pi]$ is a dense point in $\widehat{A}$, the hull of $I_n$ has empty
interior in $\widehat{A}$. Let $F_n=\{\varphi \in P(A):\varphi |I_n=0\}$.
Since the map from $P(A)$ to $\widehat{A}$ is open, we have that $F_n$ is a
closed nowhere dense set relative to $P(A)$.  It now follows from the Baire
category theorem, applied to $P(A)$, and 5.10 that there is a faithful
irreducible representation $\pi '$ whose associated pure states have 
non-empty interior in $P(A)$.  From the openness of the map from $P(A)$ to
$\widehat{A}$, we conclude that $\widehat{A}$ has an open point, whence
$A$ has an ideal $K$, necessarily essential, such that $\widehat{K}$ has
only one point.  The proof is concluded by showing $K\cong \Cal K(H)$, 
and this can be done in at least two ways.  There is a simple way
to prove that
every separable $C^*-$algebra whose spectrum is a single point must be 
elementary (i\.e\., the affirmative answer to Naimark's question in the
separable case),
or one can apply Glimm's theorem to $K$.
\enddemo

The following example demolishes one naive conjecture with regard to the
structure of type I closed faces.

\example{Example 5.12}
If $A$ is any non-type I separable $C^*-$algebra, then $A$ has a type I
closed face $F(p)$ ($p$ is even compact) such that $F(p)$ is not
isomorphic to a closed face of any type I $C^*-$algebra.

If $A$ is not unital, we consider $A^{**}$ as a subalgebra of 
$\widetilde A^{**}$ and construct $p$ as a projection in $A^{**}$ closed in
$\widetilde A^{**}$, so that $p$ will be compact.  Since $A$ is not type I,
there is an irreducible representation $\pi$ such that
$\pi(A)\not\supset\Cal K(H_\pi)$.  For the natural extension of $\pi$
to $\widetilde A$, also denoted $\pi$, we also have $\pi(\widetilde A)\not\supset
\Cal K(H_\pi)$.  Let $v_0$ be a unit vector in $H_\pi$, 
$\varphi_0=(\pi(\cdot)v_0,v_0)$ and $p_0$ the support projection in $A^{**}$
of $\varphi_0$.  By a result of Glimm \cite{{16}, Theorem 2}, there is a sequence
$\{v_n\}$ of unit vectors in $H_\pi$ such that $v_n\overset w\to\longrightarrow
0$ and $(\pi(\cdot)v_n,v_n)\to\varphi_0$ in $\widetilde A^*$.  By using the
Gram-Schmidt process, we can find a subsequence $\{v_{n_i}\}$ and an
orthonormal sequence $\{w_{n_i}\}$ such that $(w_{n_i},v_0)=0$ and
$\Vert w_{n_i}-v_{n_i}\Vert\to 0$.  Let $\varphi_i=(\pi(\cdot)w_{n_i},w_{n_i})$.

Since $p_0$ is a minimal projection in $A^{**}$, it is closed in
$\widetilde A^{**}$.  Let $B$ be the hereditary $C^*-$subalgebra of $\widetilde A$
supported by $1-p_0$, and let $e$ be a strictly positive element of $B$.
Since $\varphi_i\to\varphi_0$ in $\widetilde A^*$ and
$\varphi_{0}|_{B}=0$, $\varphi_i(e)\to 0$.  Passing to a subsequence, we may
assume $\sum \varphi_i(e)<\infty$.  Let $p_i$ be the support projection of
$\varphi_i$.  $p_i$ is in $B^{**}\cap A^{**}$, considered as a subalgebra
of $\widetilde A^{**}$.  By 0.1(ii), $\sum_1^\infty p_i$ is closed in $B^{**}$.
Thus if $p=\sum_0^\infty p_i$, $p$ is closed in $\widetilde A^{**}$.  Since
$p\in A^{**}$, $p$ is a compact projection in $A^{**}$.  Since
$pA^{**}p\cong B(H_0)$ where $H_0=\overline{\text{span}}\{v_0,w_{n_1},w_{n_2}
\dots\}$, $p$ is a type I projection.

Suppose $F(p)$ were isomorphic to a closed face, $F(p')$, of a type I
$C^*-$algebra $A'$.  Since $p'(A')^{**}p'$ can be identified with the space
of bounded affine functionals vanishing at $0$ on $F(p')$, $p'(A')^{**}p'$
is Jordan $*-$isomorphic to $pA^{**}p$.  Therefore $p'(A')^{**}p'$ is a type
I factor, and $p'$ is associated with a single irreducible representation,
$\pi'$, of $A'$.  Since $A'$ is type I, $\pi'(A')\supset \Cal K(H_{\pi'})$.
Let $\varphi_i'$, $i\ge 0$, be the element of $F(p')$ corresponding to
$\varphi_i$.  Then $\varphi_i'\to\varphi_0'$ in ${A'}^*$.  This contradicts
the facts that $\{\varphi_i'\}$ arises from an orthonormal sequence of
vectors in $H_{\pi'}$ and $\pi'(A')\supset\Cal K(H_{\pi'})$.

We think it is fairly obvious from the proof of 0.1(ii) ([12, Lemma 3]),
that the faces $F(p)$ constructed above are all isomorphic.  In Section 7
we will determine the structure of $pAp$, and this will be our formal proof
of this fact.

Finally, we want to generalize 5.5 for use in connection with a remark in
Section 7.  If $h\in A^{**}$, we say that $h$ satisfies the barycenter
formula if, when regarded as a function on $Q(A)$, $h$ is measurable with
respect to (the completion of) any regular Borel measure and
$\varphi(h)=\int h(\omega)d\mu(\omega)$ whenever $\mu$ is a regular Borel
measure on $Q(A)$ and $\varphi$ is the resultant of $\mu$.  If $A$ is
separable, it is sufficient to verify the formula for measures supported on 
$P(A)$.  Also when $A$ is separable, the barycenter formula is equivalent
to: $\pi_s^{**}(h)$ is a measurable field of operators and
$\pi^{**}(h)=\int^\oplus \pi_s^{**}(h)d\mu(s)$, whenever $\pi=\int^\oplus
\pi_s d\mu(s)$, a standard direct integral;  and again it is sufficient to
verify this in the special case where each $\pi_s$ is irreducible.  Thus
for $A$ separable the set of elements of $A^{**}$ satisfying the
barycenter formula is a $C^*-$algebra closed under weak sequential
convergence.  This $C^*-$algebra is at least as large as
$\{h:\text{Re}\ h, \text{Im}\ h \text{ are universally measurable}\}$ and appears to be
a good thing to use, though the monotone sequential closure of $A$ (discussed
in [25, \S4.5]) would do for our purposes.  For $A$ non-separable, we
know of nothing more general than the space of universally measurable
operators (\cite{26}).
\endexample

\proclaim{Theorem 5.13}
If $A$ is a separable $C^*-$algebra, $p$ is a projection in $A^{**}$
satisfying the barycenter formula, and if $\pi^{**}(p)\pi(A)\pi^{**}(p)
\subset\Cal K(H_\pi)$ for every irreducible representation $\pi$ of $A$,
then $p$ is type I.
\endproclaim

\demo{Proof}
First note that the proof of 5.5 is equally valid if $p$ satisfies the
barycenter formula instead of being closed.  Let $e$ be a strictly positive
element of $A$.  Then for $\epsilon>0$, $E_{[\epsilon,\infty)}(pep)$ satisfies
the barycenter formula and $\pi^{**}(E_{[\epsilon,\infty)}(pep))$ has finite
rank for $\pi$ irreducible.  Therefore each $p_n$ is type I where
$p_n=E_{[n^{-1},\infty)}(pep)$.  Since $p_n\nearrow p$, $p$ is type I.
\enddemo
\medskip
\noindent
{\bf 6.  More on Closed Faces with (NCEB) for $A$ Separable.}

The notations $X,p_x,\widetilde X$, and $f$ have the same meanings as in
Section 4.

\proclaim{Theorem 6.1}
If $A$ is a separable $C^*-$algebra, $p$ is a closed projection in
$A^{**}$, and if $p$ satisfies (NCEB), then $\sum_{x\in S} p_x$ is the
atomic part of a closed projection $p_S$ for every closed subset $S$
of $X$.  Also $p$ satisfies (CEB) if and only if $p_S$ is compact for $S$ compact.
\endproclaim

\demo{Proof}
By 5.6, $p$ is type I.  Let $\widetilde S=f^{-1}(S)$, a closed subset of
$\widetilde X$, let $E=\widetilde S\cap P(A)=\widetilde S\cap S(A)$, a saturated
subset of $F(p)\cap P(A)$, and let  $p_S$ be the projection called $p_E$
in 5.4.  By 5.4, $F(p_S)=\{\int \omega d\mu(\omega):\mu\text{ is a
probability measure on } E\cup\{0\}\}=$
$\{\int \omega d\mu(\omega):\mu\text{ is a
probability measure on } \widetilde S\cup\{0\}\}$.  Since $\widetilde S\cup\{0\}$
is a compact subset of $A^*$, $F(p_S)$ is closed, and hence $p_S$ is closed
by \cite{{15}, Theorem 4.8}.  By 5.4, $F(p_S)\cap P(A)=E$, and this implies that the atomic
part of $p_S$ is $\sum_{x\in S} p_x$.

If $p$ satisfies (CEB) and $S$ is compact, then $E=\widetilde S$, and $\widetilde S$
is compact by 4.3.  Thus $F(p_S)\cap S(A)=\{\int \omega d\mu(\omega):\mu
\text{ is a probability measure on } \widetilde S\}$, a closed subset of $A^*$.
Therefore $p_S$ is compact.  Conversely, if $S$ compact implies $p_S$ compact
and if $\varphi_n \to t\varphi$, for $\varphi_n$, $\varphi$ in 
$F(p) \cap P(A)$, then there is a compact set $S$ such that $f(\varphi_n) \in S$
for $n$ sufficiently large.  Since $F(p_S) \cap S(A)$ is closed, it follows
that $t=1$.
\enddemo

\example{Remarks}
1.  If $p$ satisfies only (NCEB) and $S$ is compact, then $p_S$ is nearly
relatively compact in the sense of \cite{13}.  

2.  The hypothesis of 4.4 included the assumption that $p$ be strongly
atomic, though this term was not used.  Theorem 6.1 shows that this assumption
can be dropped if $A$ is separable. 
Also the discussion after 5.7 shows that if $A$ is separable, $p$ is closed and
atomic, and $p$ satisfies (NCEB), then $p$ is strongly atomic.  Thus for
$A$ separable, the hypothesis of 4.5 can be weakened by replacing strongly
atomic with atomic.

3.  In view of the remarks after 5.4, it is not hard to calculate the
facial topology on the extreme boundary of $F(p)$, when $p$ is closed and
satisfies (CEB).  Its $T_0$-ification is the compact Hausdorff space
$X\cup\{\infty\}$.  If $p$ satisfies only (NCEB), we still see that the
closed split faces of $F(p)$ containing $0$ are in one-to-one correspondence
with the closed subsets of $X$.
\endexample
\medskip
\noindent
{\bf 7.  Some Relationships among Prior Sections and Concluding Remarks.}

Each of the three main parts of this paper (Sections 2, 3, and 4-6) studies a
different generalization of the situation considered in \cite{5} (1.4 and 7.2 below are used to justify this statement).  Sections 3-6 were motivated by our
desire to investigate the circumstances in which 2.4 applies, but the
detailed discussion below makes it clear that we have not solved this problem --
if it can be called a ``problem''.  There is a broader ``problem'' to
which all three parts of this paper are relevant:  Study the structure of
those closed faces of $C^*-$algebras which are closely modeled on locally
compact Hausdorff spaces.  We now discuss the relationships among the prior
sections.

First we consider the relationship between Sections 2 and 3.  The next result
and the remarks following show that if we were willing to use the theory of
relative $q$-continuity in the construction of MASA's, it would have been
sufficient to prove the special case of 2.4 in which the projection $p_X$
is central and abelian.  However, so far as we know, this special case is
no easier.

\proclaim{Theorem 7.1}
Let $A$ be a separable $C^*-$algebra and $p$ a closed projection in $A^{**}$.
Suppose $B$ is a commutative $C^*-$subalgebra of $SQC(p)$ which is non-degenerately embedded in $pA^{**}p$. If $\widehat B$ is totally disconnected, then there is
a commutative $C^*-$subalgebra $C$ of $A$ such that $C$ contains an
approximate identity of $A$, $p\in C^{**}$, and $pC=B$.
\endproclaim

\demo{Proof}
Let $\overline B=\{a\in A: ap=pa \text{ and } pa\in B\}$.  Then
$\text{her}(1-p)$ is an ideal of $\overline B$ and $\overline B/\text{her}(1-p)\cong B$.  We can apply 2.4 (or 2.3) with $\overline B$ playing the role of $A$ and
with $X=\widehat B$. For $x$ in $X$, $p_x$ is the support projection in
$\overline B^{**}$ of the pure state of $\overline B$ given by $x$.  Let $C$ be
the MASA of $\overline B$ given by 2.4.  Since $B$ is non-degenerate in
$pA^{**}p$, $\overline B$ hereditarily generates $A$.  Since $C$ hereditarily
generates $\overline B$ by 2.4, $C$ also hereditarily generates $A$.  One way
to deduce from 2.4 that $pC$, which is the image of $C$ in $\overline B/\text{her}
(1-p)$, is all of $B$ is to quote the classical Stone-Weierstrass theorem.
\enddemo

Suppose $p$ is a closed projection in $A^{**}$ such that $SQC(p)$ is
non-degenerate in $pA^{**}p$ (cf\. 3.2) and that $B$ is a MASA in
$SQC(p)$ which hereditarily generates $SQC(p)$.  If $A$ is separable and
$\widehat B$ is totally disconnected, then 7.1 gives a commutative algebra
$C$ (which could be assumed a MASA in $A$).  For each $x$ in $\widehat B$ we
have a pure state $\varphi_x$ of $B$ (or $C$) which is supported by a
minimal projection $p_x$ in $B^{**}$, and it follows from $pC=B$ and
$p\in  C^{**}$ that also $p_x\in C^{**}$.  If $p_x$ is minimal
in $A^{**}$, then $\varphi_x$ satisfies (UEP) relative to the inclusion of
$C$ in $A$.  If $pAp$ is an algebra (cf 3.1 and 7.2 below), then we need only
start with a MASA $B$ in $pAp$ which hereditarily generates $pAp$ and such that
each pure state of $B$ satisfies (UEP) relative to $pAp$.  It was pointed
out in Section 0 that under the hypotheses of 2.4 every element of $C_0(X)$ 
gives an element of $SQC(p_X)$.  It can be shown that $C_0(X)$ is thus
embedded as a MASA in $SQC(p_X)$ and that $C_0(X)$ is nondegenerate in
$p_XA^{**}p_X$.  Thus the above discussion applies.

\proclaim{Proposition 7.2}
Conditions (i)-(iv) of 1.2 imply that $pAp$ is an algebra, and $p$ satisfies:
\medskip
\noindent {\rm (G)}\qquad\qquad\qquad\qquad\qquad\qquad $[P(A)\cap F(p)]^-\subset [0,1
]P(A)$.
\medskip
\endproclaim

\demo{Proof}
The reduced atomic representation, $\pi$, of $A$ is faithful on
$pA^{**}p$.  Moreover, 
$\pi^{**}(pA^{**}p)\cap\Cal K(H_\pi)$ is a $C^*-$algebra, and
by 1.2(iv), $\pi^{**}(pAp)$ is contained in this algebra.  We show equality.
Let $h$ be an element of $pA^{**}_{\sa}p$ such that $\pi^{**}(h)$ is
compact.  It is sufficient to show  $h\in SQC(p)$.  If $F$ is a closed
subset of ${\Bbb R}$, then 1.2(i) implies that $E_F(h)$, computed in
$pA^{**}p$, is closed.  (In fact we don't need $F$ closed for this.)  If
$0\notin F$ (and $F$ is closed), then $\pi^{**}(E_F(h))$ is a finite rank
operator on $H_\pi$, and by \cite{{1}, Corollary II.8} this implies $E_F(h)$ is compact.  Thus
$h\in SQC(p)$.  Then (G) follows from [17, \S5] for example.
\enddemo

The same reasoning shows $pA^{**}p=QC(p)$.

Next we consider the relationship between Sections 2 and 4.  By 6.1 if $A$
is separable and $p$ is a closed projection in $A^{**}$ satisfying (CEB)
then we have the hypotheses of 2.3 (except for total disconnectedness of
$X$).  By 4.1 each $p_x$ is of finite rank.  For 2.4, we would want each
$p_x$ to be of rank 1.  This happens for the $p_x$'s of Section 4 if and
only if $p$ is abelian.  If $p$ is not abelian, it might be possible to
write $p_x=p_{x,1}+\dots+p_{x,n_x}$ where the $p_{x,i}$'s are minimal and
$\{p_{x,i}\}$ satisfies the hypotheses of 2.4 with $X$ replaced by some
space $\overline X$.  Then $\overline X$ would map onto $X$ by a closed continuous
map with finite fibers.  However, Example 7.6(a) below shows that this is not
always possible.

Conversely, suppose the hypotheses of 2.4 are satisfied.  By 1.4, if $X$ is
countable and discrete and each equivalence class of $\{\varphi_x:x\in X\}$ is finite, then
$p_X$ satisfies (CEB).  If $p_X$ is abelian, we can deduce (7.3 below) that
$p_X$ satisfies (CEB) even for $X$ not discrete;  but it is fairly obvious
(cf 7.6(b) below) that in general $p_X$ need not satisfy
(NCEB) or even (G).

In retrospect it seems that (G) is worthy of more study in the present context
despite the fact, as pointed out in the remark following 4.5, that it does
not imply a positive answer to the isolated point question.  One reason is
mentioned in 7.2 above.
However, the conclusion of 1.4 is definitely false if
we drop the hypothesis of finite equivalence classes.  (This follows from 4.1.)
It may be that (G) is part of the hypothesis of a nice result.  Also, even
though, by Example 3.4, (G) does not imply that $F(p)$ is associated with a
locally compact Hausdorff space, we do not know whether (G) implies that
$F(p)$ is associated with a Hausdorff space.  We will show below that (G) does
imply that $p$ is type I.

One could also consider weaker conditions than (G):
$$
[P(A)\cap F(p)]^-\subset\{\text{type I factorial quasi-states}\}\tag1
$$
$$
[P(A)\cap F(p)]^-\subset z_{\at}A^*.\tag2
$$
(2) is suggested by the theory of perfect $C^*-$algebras (\cite{8}).

With regard to the relationship between Sections 3 and 4, we note that a
closed projection $p$ satisfying (CEB) need not satisfy (MSQC) 
(cf. 7.6(c) below).  However, it follows from 6.1 that $p$ does satisfy
the hypothesis of 3.2 (A separable).  Also if $p$ is closed and abelian
and satisfies (CEB), then $p$ does satisfy (MSQC).  (It follows from 5.3
that $F(p)$ is isomorphic to the set of probability measures on
$\widetilde X\cup\{0\}$.  Is there a less technical proof?)  It may be that
there are other useful concepts on the extensiveness of $SQC(p)$.

At the end of this section we return to Example 5.12, partly to show that
it does satisfy the conclusion of 4.5.  A complete theory for closed faces
of $C^*-$algebras analogous to the theory of scattered $C^*-$algebras might 
have to be quite complicated.

\proclaim{Proposition 7.3}
Assume the hypotheses of 2.4 and also that $p_X$ is abelian.  Then $p_X$
satisfies (CEB).
\endproclaim

\example{Remark}
The hypothesis that $p_X$ is abelian can be rephrased more concretely:  The
$\varphi_x$'s are mutually inequivalent. (\cite{10}).
\endexample

\demo{Proof}
Since $p_X$ is abelian, $P(A)\cap F(p_X)=\{\varphi_x:x\in X\}$.  Suppose
$\varphi_{x_i}\to\psi$ in $Q(A)$.  Passing to a subnet, we may assume
$x_i\to x$ in $X$ or $x_i\to\infty$.  If $x_i\to x$, let $\{S_j\}$ be a set
of compact neighborhoods of $x$ such that $\bigcap_j\ S_j=\{x\}$.
Since $p_{S_j}$ is compact and $\varphi_{x_i}\in F(P_{S_j})$ for $i$ large,
we conclude that $\Vert\psi\Vert=1$ and $\psi\in \bigcap_j\ F(p_{S_j})
=F(\bigwedge_j\ p_{S_j})$.  Since a closed projection is
determined by its atomic part, $\bigwedge_j\ p_{S_j}=p_x$, and hence
$\psi=\varphi_x$.  If $x_i\to\infty$, let $\{U_j\}$ be a set of relatively
compact open subsets of $X$ such that $\bigcup_j\ U_j=X$, and
let $S_j=X\setminus U_j$.  Since $x_i$ is eventually in $S_j$ and since
$p_{S_j}$ is closed, we see that $\psi\in \bigcap_j\ F(P_{S_j})=
F(\bigwedge_j\ p_{S_j})=\{0\}$.
\enddemo

It would be desirable if the hypotheses in 2.4 that certain projections
are atomic parts of closed projections could be
stated entirely in terms of pure
states (or equivalently, minimal projections).
This can be done in a situation of intermediate generality. Consider the
following conditions for a projection $p$ in $A^{**}$:
$$
\exists t\in (0,1]\ \text{such that}\ [P(A)\cap F(p)]^-\subset \{0\} \cup [t,1][P(A)\cap F(p)].\tag3
$$
$$
\{0\} \cup [P(A)\cap F(p)] \ \text{is closed.}\tag4
$$
$$
[P(A)\cap F(p)] \ \text{is closed.}\tag5
$$
Conditions (3) and (4) are strengthenings of (NCEB) and (CEB) respectively,
and are equivalent to (NCEB) and (CEB) if $p$ is closed.  But we are
interested in the case where $p$ is atomic.  If $p$ is the atomic part of a
closed projection $q$, then $q$ satisfies (NCEB) or (CEB) if and only if $p$
satisfies (3) or (4).  If $p$ is atomic and satisfies (3) or (4), is $p$
necessarily the atomic part of a closed projection?  We can prove this if $p$
is strongly atomic, in which case $p$ itself is closed.  (In general let
$C$ be the closed convex hull of $\{0\} \cup [P(A)\cap F(p)]$.  If $q$ exists,
then $F(q)=C$.  The tricky thing is to prove that $C$ is a face of $Q(A)$.)

\proclaim{Lemma 7.4}
If $p$ is a strongly atomic projection satisfying (3), then $p$ is closed.
\endproclaim

\demo{Proof}
Let $X=\{0\} \cup [P(A)\cap F(p)]^-$ and let $C$ be the closed convex hull
of $X$.  Then every element of $C$ is the resultant of a probability measure
on $X$.  Since $X$ is norm separable, the resultant is actually a Bochner
integral.  Hence $C \subset F(p)$.  The reverse inclusion follows easily from
the structure of atomic von Neumann algebras.
\enddemo

\example{Remark}
The same argument works if (3) is replaced by a similar modification of (2).
\endexample

\proclaim{Corollary 7.5}
Let $X$ be a locally compact Hausdorff space with only countably many points,
and let $\{p_x: x \in X\}$ be a family of mutually orthogonal minimal
projections in $A^{**}$. If $\sum_{x\in S} p_x$ satisfies (3) whenever $S$
is a closed subset of $X$ and (5) when $S$ is compact,
and if $A$ is separable, then we have the
hypotheses and conclusions of Corollary 2.4.
\endproclaim

\example{Remarks}
1.  Let $\varphi_x$ be the pure state supported by $p_x$.  If $p=\sum_{x \in X} p_x$
is abelian, i.e. if the $\varphi_x$'s are mutually inequivalent, then the
hypotheses on $\sum_{x\in S} p_x$ can be stated more concretely:
$$
\text{If}\  x_n \to x,\ \text{then}\  \varphi_{x_n} \to \varphi_x,\ \text{and if}\  x_n \to \infty,\ \text{then}\  \varphi_{x_n} \to 0.\tag6
$$
We can actually replace the hypotheses on $\sum_{x\in S} p_x$ by (6) if we assume
only that the equivalence classes have bounded finite cardinality.  (Note that
they have to be finite by 4.1 if $p$ satisfies (3).)  The proof of this uses
Akemann's result in \cite{1} that the supremum of finitely many mutually 
commuting closed projections is closed.

2.  Even when $X$ is uncountable, the hypotheses on $\sum_{x \in S} p_x$ in
2.4 can be modified somewhat:  If $\sum_{x \in X} p_x$
is the atomic part of a closed projection, and if $\sum_{x \in S} p_x$ satisfies
(3) for $S$ closed and (5) for $S$ compact, then we have the hypotheses of 2.4.
We will not provide a complete proof of this because it would be rather 
technical and it is not clear that the result is a big improvement on 2.4.
The main lemma is the following:  

\noindent Let $p$ be a closed projection satisfying (NCEB) and $q$ a
subprojection of $z_{\at}p$.  If $A$ is
separable and $q$ satisfies (3), then $q$ is the atomic part
of a closed projection.

\noindent The proof of this uses 5.6, the other results of Section 5 (in 
particular the discussion following 5.8), and the von Neumann selection lemma.
\endexample

\example{Examples 7.6}
(a) Let $A=c\otimes\Cal K$ and define a closed projection $p$ in
$A^{**}$ by letting $p_\infty$ be the projection on $\text{span}\{e_1,e_2\}$ and
$p_n$ the projection on
$$
\cases  {\Bbb C}e_1, &\quad n=3k+1\\
{\Bbb C}e_2, &\quad n=3k+2\\
{\Bbb C}(2^{-\frac12}e_1+2^{-\frac12}e_2), &\quad n=3k\endcases
$$
It is easy to see that $p$ satisfies (CEB).  Let $\varphi_n$ be the pure state
of $A$ supported by $p_n$, $n<\infty$, and suppose $B$ is a MASA of $A$ such
that each $\varphi_n|_B$ satisfies (UEP).  Thus each $p_n$ is in $B^{**}$.
If $b\in B$, then $e_1$ is an eigenvector of each $b_{3k+1}$ and hence
$e_1$ is an eigenvector of $b_\infty$.  Similarly $e_2$ and $2^{-\frac12}e_1
+2^{-\frac12}e_2$ are eigenvectors of $b_\infty$.  Therefore all three
eigenvalues are the same and $b_\infty p_\infty=\lambda p_\infty$.  It
follows that $p_\infty$ is a minimal projection of $B^{**}$.  Thus no
matter how we write $p_\infty=p'+p''$, with $p'$ and $p''$ rank one
projections, we cannot achieve the conclusion of 2.4, let alone the hypotheses.

(b)  First note that if $p$ is the projection of 5.12, then we have the hypotheses
of 2.4 with $X={\Bbb N}\cup\{\infty\}$ and $p=p_X$.  Since all of the
$\varphi_n$'s, $1\le n\le\infty$, are equivalent, it is easy to see that the
non-pure state $\tfrac 12 \varphi_1+\tfrac 12 \varphi_\infty$ is in 
$[P(A)\cap F(p)]^-$ (cf \cite{16}), so that $p$ does not satisfy (G).

It is better to give an example where the equivalence classes of
$\{\varphi_x:x\in X\}$ are finite, since by 4.1 there is no hope of (NCEB)
without this finiteness.  A standard example suffices for this.  Let
$A=\{a\in c\otimes M_2:a_\infty \ \text{ is diagonal}\}$.   Let
$B=\{a\in A: a_n \text{ is diagonal, } \forall n\}$. 
Then $B$ is a
MASA in $A$, and we let $X=\widehat B$, the disjoint union of two copies of
${\Bbb N}\cup\{\infty\}$.  It is clear that for $x$ in $X$ the pure state
$\psi_x$ of $B$ satisfies (UEP);  and if $p_x$ is the support projection
of $\psi_x$, we have the hypothesis of 2.4 with $p_X=1_A$.  Since $\widehat A$
is not Hausdorff, it follows from [17, Thm. 6] that $p_X$ does not satisfy (G).
Of course, this is also easy to see explicitly.

It is possible to give a similar example in which $\widehat A$ is Hausdorff, but a
different condition of [17, Theorem 6] is violated.  Let
$A=\{a\in c\otimes M_2\otimes M_2:a_\infty\in M_2\otimes I_2\}$.  If
$B=\{a\in A:a_n\in D_2\otimes D_2,n<\infty; a_\infty\in D_2\otimes I_2\}$, where
$D_2=\{d\in M_2:d \text{ is diagonal}\}$, then $B$ is a MASA in $A$ and we
can proceed similarly to the above.  Again $X$ is the disjoint union of two
copies of ${\Bbb N}\cup\{\infty\}$ (arising more naturally as
${\Bbb N}\cup{\Bbb N}\cup\{\infty\}$).

(c)  Consider one of the examples of alternative 2 of 4.6 constructed in 4.8(c).
Let $k=n=2$ and $T=1$.  We then get a projection $p$ satisfying (CEB) where
the space $X$ of Section 4 is ${\Bbb N}\cup\{\infty\}$, rank $p_n=2$,
$1\le n\le\infty$, and there is a single element $\varphi$ of
$F(p_\infty)\cap P(A)$ such that every sequence $(\psi_n)$ with $\psi_n$
in $F(p_n)\cap P(A)$ converges to $\varphi$.  In this case we can write
$p_n=p_{n,1}+p_{n,2}$ so that the hypotheses of 2.4 are satisfied.  All we have
to do is take $p_{\infty,1}$ to be the support projection of $\varphi$.
$\widetilde X$ will be homeomorphic to ${\Bbb N}\cup\{\infty\}$, but it arises as 
the disjoint union of ${\Bbb N}\cup {\Bbb N}\cup\{\infty\}$
with an isolated point.  This
example does not satisfy (MSQC).  One way to see this is to note that the
saturation of an open subset of $F(p)\cap P(A)$ need not be open, and hence
$F(p)$ is not isomorphic to the quasi-state space of a $C^*-$algebra.
Explicitly, any element $h$ of $SQC(p)$ (or $QC(p)$) must have $\varphi$
definite on $h_\infty$.
\endexample

\proclaim{Lemma 7.7}
If $A$ is a $C^*-$algebra, $p$ is a projection in $A^{**}$, and $p$ satisfies
(G), then $\pi^{**}(pAp)\subset\Cal K(H_\pi)$ for every irreducible
representation $\pi$ of $A$.
\endproclaim

\demo{Proof}
Part of the proof of 4.1 applies:  If $(e_n)$ is an orthonormal sequence
in the range of $\pi^{**}(p)$ and $\psi_n=(\pi(\cdot)e_n,e_n)$, we can
conclude that $\psi_n\to 0$.  If $E_{[\epsilon,\infty)}(\pi^{**}(pap))$ has
infinite rank for some $a$ in $A_+$ and $\epsilon>0$, then we can obtain a
contradiction by taking the $e_n$'s in the range of $E_{[\epsilon,\infty)}
(\pi^{**}(pap))$.
\enddemo

\proclaim{Corollary 7.8}
If $A$ is a separable $C^*-$algebra, $p$ is a projection in $A^{**}$
satisfying the barycenter formula (in particular if $p$ is closed), and $p$
satisfies (G), then $p$ is type I.
\endproclaim

\demo{Proof}
Apply 5.13.
\enddemo

\example{7.9. Continuation of Example 5.12}

(a)  A subprojection $p'$ of $p$ is closed if and only if $p'$ has finite
rank or $p'\ge p_0$.

{\it Proof}.
If $p'$ has finite rank, then $p'$ is closed by \cite{1}.  If $p'\ge p_0$, then $p'-p_0$ is closed in $B^{**}$ by 0.1(ii).  Therefore $p'$ is closed in
$\widetilde A^{**}$ and {\it a fortiori} in $A^{**}$.

If $p'$ has infinite rank, then the range of $p'$ contains an infinite
dimensional subspace $H'$ of $\overline{\text{span}}\{w_{n_1},w_{n_2},\dots\}
=\text{ range}(p-p_0)$.  (This last is a codimension 1 subspace of
$H_0=\text{ range }p$.)  Let $(u_n)$ be a sequence of unit vectors in
$H'$ such that $u_n\overset w\to\longrightarrow 0$ and 
$\psi_n=(\pi(\cdot)u_n,u_n)$. Then $\psi_n|_B\to 0$, since $\pi^{**}(pBp)
\subset\Cal K(H_\pi)$.  Since $p$ is compact, it follows that
$\psi_n\to\varphi_0$.  If $p'$ is closed, this implies $\varphi_0\in F(p')$
and hence $p'\ge p_0$.

(b) For any non-zero closed subprojection $p'$ of $p$ there is a minimal
projection $p_1$ such that $p_1\le p'$ and $p'-p_1$ is closed.

{\it Proof}.
If $p'$ has infinite rank, then $p'\ge p_0$.  We can find a minimal projection
$p_1$ such that $p_1\le p'-p_0$.  Then $p_0\le p'-p_1$ so that $p'-p_1$
is closed.  If $p'$ has finite rank then $p'-p_1$ is closed for any
choice of $p_1$.

(c)  If $pA^{**}p$ is identified with $B(H_0)$, then 
$$
\aligned
pAp &=\{x\in B(H_0):x-\varphi_0(x)I_{H_0}\in\Cal K(H_0)\}\\
&=\{x\in B(H_0):x-(xv_0,v_0)I_{H_0}\in\Cal K(H_0)\}.
\endaligned
$$

{\it Proof}.
Let $H_1=H_0\ominus{\Bbb C}v_0$.  By construction and the proof of 7.2,
applied to $B$ and $p-p_0$,
$pBp=\Cal K(H_1)$.  Note that since $p$ is compact, $pAp=p\widetilde A p$.
To show that $pAp$ is contained in the set indicated, it is enough to
show $pap$ is compact when $a\in\widetilde A$ and $\varphi_0(a)=0$.
By [25, 3.13.6], $a=l+r$, where $l\in\widetilde AB$ and $r\in B\widetilde A$.
Since $x$ is compact if and only if $x^*x$ is compact, $Bp\subset\Cal K(H_\pi)$;  and similarly $pB\subset\Cal K(H_\pi)$.  Therefore $pap\in\Cal K(H_0)$.

For the reverse inclusion, since $p\in pAp$, $\Cal K(H_1)\subset pAp$, and
$pAp$ is self-adjoint, it is sufficient to show that $pAp$ contains every
rank 1 operator $x$ of the form $v\to (v,v_0)v_1$ for $v_1\in H_1$.  By the
Kadison transitivity theorem (\cite {20}) there is $a$ in $A$ such that $av_0=v_1$ and
$a^*v_0=0$.  By the above, since $(av_0,v_0)=0$, $pap$ is compact.  Hence
$(pap-x)\in\Cal K(H_1)=pBp$.  This implies that $x$ is in $p\widetilde A p$.

Since by [7, 4.4], the bidual of $pAp$ is $pA^{**}p$, and since the predual of
a $W^*-$algebra is unique, it follows from (c) that the dual space of $pAp$
is $\Cal T(H_0)$, the set of trace class operators on $H_0$.  A concrete
statement of this reads:

If $v_0$ is a unit vector in the (separable, infinite dimensional) Hilbert
space $H_0$ and 
$$
A_{v_0}=\{x\in B(H_0):x-(xv_0,v_0)I_{H_0}\ \text{ is compact}\},
$$
then the dual space of $A_{v_0}$ is naturally isometrically isomorphic to
$\Cal T(H_0)$.  In particular, for $T\in\Cal T(H_0)$,
$\Vert T\Vert_1=\sup\{|Tr(Tx)|:x\in A_{v_0},\Vert x\Vert\le 1\}$.

It is amusing to give a direct proof of this statement (which removes the
parenthetical part of the hypothesis).  The main step is to prove the second
sentence.

Finally, we note that 5.12 gives another example of how the behavior of
closed faces of $C^*-$algebras differs from that of $C^*-$algebras.  If
$\pi$ is an irreducible representation of a $C^*-$algebra $A$, then
$\pi(A)\cap\Cal K(H_\pi)$ is either $0$ or $\Cal K(H_\pi)$.  The
analogous statement for Example 5.12 (replacing $A$ by $pAp$) is false.
\endexample

\vfill\eject

\vskip.25truein
\hskip3truein Dept. of Mathematics,

\hskip3truein Purdue University

\hskip3truein W. Lafayette, IN 47907, USA

\hskip3truein e-mail: lgb\@math.purdue.edu.

\enddocument